\pdfoutput=1
\RequirePackage{ifpdf}
\ifpdf
\documentclass[pdftex]{sigma}
\else
\documentclass{sigma}
\fi

\numberwithin{equation}{section}

\newtheorem{Theorem}{Theorem}[section]

\newtheorem{Lemma}[Theorem]{Lemma}
\newtheorem{Proposition}[Theorem]{Proposition}
{
 \theoremstyle{definition}
 \newtheorem{Definition}[Theorem]{Definition}
 
 \newtheorem{Remark}[Theorem]{Remark}
 \newtheorem{Notation}[Theorem]{Notation} 
 \newtheorem{constructionPrinciple}[Theorem]{Construction Principle} 
}

\usepackage{amsxtra}
\usepackage{amscd}
\usepackage{tikz}
\usepackage{tikz-cd}

\newcommand{\df}{\it }

\DeclareMathOperator{\cocurvature}{cocurv}

\DeclareMathOperator{\SO}{\operatorname{SO}}

\DeclareMathOperator{\rank}{rank} 

\DeclareMathOperator{\adjoint}{ad}
\DeclareMathOperator{\Adjoint}{Ad}

\DeclareMathOperator{\trace}{trace}

\DeclareMathOperator{\dimension}{dim}

\DeclareMathOperator{\automorphism}{Aut}

\DeclareMathOperator{\spann}{span}
\DeclareMathOperator{\curvature}{curv}
\DeclareMathOperator{\torsion}{tor}
\DeclareMathOperator{\rad}{rad}

\newcommand{\ii}{\mathrm{I\!I}}

\newcommand\monomorphism{\hookrightarrow}

\newcommand\suchthat{\, \vert \,}

\renewcommand{\ge}{\geqslant}
\newcommand*\circled[1]{{\tiny\tikz[baseline=(char.base)]{
 \node[shape=circle,draw,inner sep=2pt] (char) {#1};}}}

\newcommand{\kk}{\varkappa}

\begin{document}
\allowdisplaybreaks

\newcommand{\arXivNumber}{1703.03851}

\renewcommand{\PaperNumber}{062}

\FirstPageHeading

\ShortArticleName{Lie Algebroid Invariants for Subgeometry}

\ArticleName{Lie Algebroid Invariants for Subgeometry}

\Author{Anthony D.~BLAOM}

\AuthorNameForHeading{A.D.~Blaom}

\Address{Waiheke Island, New Zealand}
\Email{\href{mailto:anthony.blaom@gmail.com}{anthony.blaom@gmail.com}}
\URLaddress{\url{https://ablaom.github.io}}

\ArticleDates{Received November 15, 2017, in final form June 13, 2018; Published online June 18, 2018}

\Abstract{We investigate the infinitesimal invariants of an immersed submanifold~$\Sigma $ of a~Klein geometry $M\cong G/H$, and in particular an invariant filtration of Lie algebroids over~$\Sigma $. The invariants are derived from the logarithmic derivative of the immersion of~$\Sigma $ into~$M$, a complete invariant introduced in the companion article, {\em A~characterization of smooth maps into a homogeneous space}. Applications of the Lie algebroid approach to subgeometry include a new interpretation of Cartan's method of moving frames and a novel proof of the fundamental theorem of hypersurfaces in Euclidean, elliptic and hyperbolic geometry.}

\Keywords{subgeometry; Lie algebroids; Cartan geometry; Klein geometry; differential invariants}

\Classification{53C99; 22A99; 53D17}

\section{Introduction}\label{introduction}%
We initiate a general analysis of subgeometry invariants using the language of Lie algebroids, which is well suited to the purpose, according to \cite{Blaom_F}. By a {\df subgeometry} we mean an immersed submanifold $\Sigma $ of a {\df Klein geometry}, which shall refer to any smooth manifold $M$ on which a Lie group $G$ is acting smoothly and transitively. The geometry will be required to be {\df simple} in the sense that its isotropy groups are weakly connected (see Definition \ref{weakD} below). For example, the Riemannian space forms ${\mathbb R}^n $, $S^n$ and ${\mathbb H}^n$ are all simple Klein geometries if we take $G$ to be the group of orientation-preserving isometries, but even if $G$ is the group of {\em all} isometries.

Informally, an object associated with $\Sigma $ is an {\df invariant} if it is unchanged when $\Sigma $ is replaced by its image under a~symmetry of~$M$. The group of symmetries of $M$ is generally larger than the transformations defined by $G$ (see Definition~\ref{symmetriesD}). However, in most cases of interest to geometers this transformation group is extended by a factor of, at most, finite order.

\subsection*{The logarithmic derivative of an immersion}
According to \cite{Blaom_F}, a complete infinitesimal invariant of a smooth map $f \colon \Sigma \rightarrow M$ is its {\df logarithmic derivative} $\delta f \colon A(f) \rightarrow {\mathfrak g} $. Here ${\mathfrak g} $ is the Lie algebra of $G$, $A(f)$ is the pullback (in the category of Lie algebroids) of the action algebroid ${\mathfrak g} \times M$ under $f$, and $\delta f $ is the composite of the natural map $A(f) \rightarrow {\mathfrak g} \times M$ with the projection ${\mathfrak g} \times M \rightarrow {\mathfrak g} $, both of which are Lie algebroid morphisms. In saying that $\delta f$ is {\df complete}, we mean that $\delta f $ determines $f$ up to symmetry.

Specialising to the case of an immersion $f \colon \Sigma \rightarrow M $, we may describe $A(f)$ more concretely as follows:
\begin{gather}\label{parlour}
 A(f)=\big\{(\xi,x)\in {\mathfrak g} \times \Sigma \suchthat \xi^\dagger (f(x)) \in T_{f(x)}\Sigma\big\},
\end{gather}
where $\xi^\dagger$ denotes the infinitesimal generator of $\xi \in {\mathfrak g} $. That is, $A(f)$ is a subbundle of the trivial vector bundle ${\mathfrak g} \times \Sigma $ over $\Sigma $ encoding which infinitesimal generators of the action of $G$ on the ambient manifold $M$ are tangent to the immersed submanifold $f(\Sigma)$ at each point. The logarithmic derivative of $f$ is just the composite $\delta f \colon A(f) \hookrightarrow {\mathfrak g} \times \Sigma \rightarrow {\mathfrak g} $.

By construction, we obtain a natural vector bundle morphism $\# \colon A(f) \rightarrow T \Sigma $, the {\df anchor} of~$A(f)$, whose kernel ${\mathfrak h} $ is a Lie algebra bundle. What may be less obvious to the general reader is that the corresponding bracket on {\em sections} of~${\mathfrak h} $ extends to a bracket on sections of~$A(f)$ generalizing the Leibniz identity for vector field brackets:
\begin{gather*}
 [X, fY] = f[X, Y] + {\rm d}f(\#X)Y, \qquad X,Y \in \Gamma(A(f)).
\end{gather*}
In other words, $A(f)$ is a Lie algebroid. This bracket is given by
the formula
\begin{gather*}
 [X,Y]=\nabla_{\#X}Y-\nabla_{\#Y}X +\{X,Y\},
\end{gather*}
which here plays a role similar to that of the classical Maurer--Cartan equations in Cartan's approach to subgeometry \cite{Blaom_F}. Here $\nabla $ is the canonical flat connection on ${\mathfrak g} \times \Sigma$ and, viewing sections of ${\mathfrak g} \times \Sigma $ as ${\mathfrak g} $-valued functions, $\{X,Y\}(x):=[X(x),Y(x)]_{\mathfrak g} $.

The logarithmic derivative of an arbitrary smooth map $f \colon \Sigma \rightarrow M$ generalizes \'Elie Cartan's logarithmic derivative of a smooth map $f \colon \Sigma \rightarrow G$ into a Lie group, and the theory laid down here is closely related to Cartan's method of moving frames, as explained in the introduction to~\cite{Blaom_F}. Indeed, Cartan's method, which cannot always be applied globally, can be reinterpreted within the new framework (Section~\ref{appendixA}) but the new theory can also be applied globally and without fixing coordinates or frames, as we shall demonstrate.

\subsection*{Aims and prerequisites}
The main contribution of the present article is to `deconstruct' logarithmic derivatives sufficiently that known invariants may be recovered in some familiar but non-trivial examples, and to show how subgeometries in these examples may be reconstructed from their invariants by applying the general theory. The analysis of the finer structure of logarithmic derivatives given here is by no means exhaustive.

Specifically, we shall recover the well-known fundamental theorem of hypersurfaces, or {\df Bonnet theorem}, for Euclidean, elliptic and hyperbolic geometry. In principle, Bonnet-type theorems for other Klein geometries could be tackled using the present framework, although this is not attempted here.

Introductions to the theory of Lie groupoids and Lie algebroids are to be found in \cite{CannasdaSilva_Weinstein_99,Crainic_Fernandes_11,Dufour_Nguyen_05,Mackenzie_05}. In particular, the reader should be acquainted with the {\df representations} of a Lie algebroid or groupoid, which will be ubiquitous. Some of the ideas presented here and in~\cite{Blaom_F} are also sketched in \cite{Blaom_H}, which may be regarded as an invitation to Lie algebroids for the geometer who is unfamiliar with (or skeptical of) these objects.

The main result needed from our companion article is a specialisation to immersions of the infinitesimal characterization of arbitrary smooth maps into a simple Klein geometry \cite[Theorem~2.16]{Blaom_F}, which we state here as Theorem \ref{infirmT}. While not essential, we recommend readers acquaint themselves with the first two sections of~\cite{Blaom_F}. Details regarding monodromy may be skimmed, as it is mostly ignored here, but closer attention should be paid to the notion of symmetry, which is subtle and not the usual one.

\subsection*{Invariants}
We now briefly describe the two types of invariants to be derived from the logarithmic derivative of an immersion $f \colon \Sigma \rightarrow M$, when $M $ is a simple Klein geometry. Details are given in Section~\ref{invint}. The first kind of invariant is an ${\mathbb R} $-valued (or ${\mathbb C}$-valued) function on $A(f)$ obtained, very simply, by composing the logarithmic derivative $\delta f \colon A(f) \rightarrow {\mathfrak g} $ with a polynomial on the Lie algebra ${\mathfrak g} $ invariant under the adjoint action. Many familiar invariants of submanifold geometry can be derived with the help of such invariants.

To describe the second kind of invariant, identify $A(f)$ with a subbundle of $E={\mathfrak g} \times \Sigma$, as in~\eqref{parlour} above. Then, unless $\Sigma $ enjoys a lot of symmetry (see below) the canonical flat connection $\nabla $ on $E$ does not restrict to a~connection on~$A(f)$. However, under a suitable regularity assumption detailed in Section \ref{invint}, there will be a largest subbundle $A(f)^2 \subset A(f)$ such that for any vector field $U$ on $\Sigma$, $\nabla_U X$ will be a section of~$A(f)$ whenever $X$ is a section of~$A(f)^2$. The new bundle~$A(f)^2$ need not be $\nabla $-invariant either, so we repeat the process and obtain an invariant filtration
\[A(f)\supset A(f)^2 \supset A(f)^3\supset \cdots \]
by what turn out to be subalgebroids of $A(f)$.

By definition, $A(f)$ is the Lie algebroid of the pullback ${\mathcal G}(f)$ of the action groupoid $G \times M$ to $\Sigma $:
\begin{gather}
 {\mathcal G}(f) = \{(y,g,x) \in \Sigma \times G \times \Sigma \suchthat f(y)=g \cdot f(x)\}.\label{ssaa}
\end{gather}
Tangent-lifting the action of $G$ on $M$, we obtain an action of the groupoid ${\mathcal G}(f)$ on $T_\Sigma M$, the vector bundle pullback of $TM$. If we define $ {\mathcal G}(f)^2 \subset {\mathcal G}(f)$ to be the isotropy of $T\Sigma\subset T_\Sigma M$ with respect to this representation, then its Lie algebroid is in fact $A(f)^2$ (Proposition~\ref{nungP}). In particular, $A(f)^2$ is the isotropy subalgebroid under the corresponding infinitesimal representation.

In fact, all the Lie algebroids $A(f)^j$ in the filtration that are transitive can be characterized as isotropy subalgebroids, under representations built out of one {\df fundamental representation} of~$A(f)$ on $E$ that we define in Section~\ref{invint}. Corresponding interpretations of the globalisations~${\mathcal G}(f)^j$, for $j \ge 3$, are not offered here, but nor are they needed for making computations.

There are evidently many ways to combine polynomial invariants with the invariant filtration to construct new invariants. For example, if $A(f)^k \cong T \Sigma $ for sufficiently large $k$, then a~polynomial invariant restricts to an ordinary symmetric tensor on $\Sigma $ that must also be an invariant. The curvature of a curve in the Euclidean, hyperbolic, elliptic or equi-affine plane may be understood in this way. For codimension-one submanifolds in higher dimensional elliptic or hyperbolic geometry (see Section~\ref{bonnetsec}) the Killing form on ${\mathfrak g} $ encodes the first fundamental form (inherited metric) of an immersion $f \colon \Sigma \rightarrow M$, and also allows one to define a normal to $A(f)$ in $E={\mathfrak g} \times \Sigma $ whose ordinary derivative, as a~${\mathfrak g} $-valued function, determines the second fundamental form of the immersion.\footnote{The normal can also be constructed {\em up to scale} by considering certain representations associated with the invariant filtration, i.e., without consideration of polynomial invariants.} In Euclidean geometry, summarized below, a similar role is played by an appropriate $\Adjoint$-invariant quadratic form on $\rad {\mathfrak g} \subset {\mathfrak g} $, the subalgebra of constant vector fields.

Finally, for each $j \ge 2$, the Lie algebroid $A(f)^j$ in the filtration acts on the tangent bundle~$T \Sigma $, which may give $\Sigma $ an intrinsic `infinitesimal geometric structure' \cite{Blaom_12}. In the Riemannian geometries mentioned above the inherited metric (up to scale) can be alternatively understood in this way.

\subsection*{Primitives and their existence and uniqueness}
Abstracting the properties of the logarithmic derivatives of smooth maps $f \colon \Sigma \rightarrow M$, we obtain generalized Maurer--Cartan forms \cite{Blaom_F}. If we restrict attention to the logarithmic derivatives of {\em immersions}, then as a special case we obtain {\df infinitesimal immersions}, formally defined in Section~\ref{infch}. An infinitesimal immersion is a Lie algebroid morphism $\omega \colon A \rightarrow {\mathfrak g} $, where $A$ is a~Lie algebroid over $\Sigma$, which is, in particular, injective on fibres. The construction of polynomial invariants and the invariant filtration given above generalizes immediately to infinitesimal immersions.

A smooth map $f \colon \Sigma \rightarrow M$ is a {\df primitive} of an infinitesimal immersion $\omega \colon A \rightarrow {\mathfrak g} $ if $\omega $ `is' the logarithmic derivative of $f$ in the following sense: there exists a~Lie algebroid morphism $\lambda \colon A \rightarrow A(f)$ covering the identity on $\Sigma $, and an element $l \in G$, such that the following diagram commutes:
\begin{gather*}
 \begin{CD}
 A @>{\omega}>> {\mathfrak g} \\
 @V{\lambda}VV @VV{\Adjoint_l}V \\
 A(f) @>{\delta f}>> {\mathfrak g}.
 \end{CD}
\end{gather*}
The primitive is {\df principal} if we can take $l=1_G$. Primitives of infinitesimal immersions are necessarily immersions, and the morphism $\lambda $ automatically an isomorphism. Specialising \cite[Theorem~2.16]{Blaom_F} to infinitesimal immersions, we have:
\begin{Theorem}\label{infirmT} Let $M$ be a Klein geometry with transitively acting group $G$. Then an infinitesimal immersion $\omega \colon A \rightarrow {\mathfrak g} $ admits a primitive $f \colon \Sigma \rightarrow M$ if and only if it has trivial monodromy, in which case $A $ is an integrable Lie algebroid. Assuming $M$ is a simple Klein geometry, the primitive $f$ is unique up to symmetry, and can be chosen to be principal.
\end{Theorem}

As explained at the beginning of this section, the meaning of `symmetry' is not quite the usual one. The {\df monodromy} is a map $\pi_1(\Sigma, x_0) \rightarrow M$, whose definition is given in \cite{Blaom_F}. The monodromy is {\df trivial} if it is constant, which is automatically true if $\Sigma $ is simply-connected (always the case in the illustrations to be given here).

\subsection*{Riemannian subgeometry and the Bonnet theorem}
As a concrete illustration of the general theory to be developed, we consider, in Sections~\ref{bonnetsec2} and~\ref{bonnetsec}, the case that $M$ is Euclidean space ${\mathbb R}^n $, the sphere $S^n$, or hyperbolic space ${\mathbb H}^n$. The problem is to classify all oriented codimension-one immersed submanifolds $\Sigma \subset M$, up to actions of the group $G$ of orientation-preserving isometries. We identify the Lie algebra ${\mathfrak g} $ of $G$ with the space of Killing fields on $M$. Our purpose is not to demonstrate anything new about Riemannian geometry, but to show how to apply and interpret the general theory in familiar cases. The following is a synopsis of the case $M = {\mathbb R}^n$.

It may be observed that the pullback ${\mathcal G}(f)$ of the action groupoid $G \times M$ by an immersion $f \colon \Sigma \rightarrow M$ acts faithfully on $T_\Sigma M$~-- that is, a rigid motion mapping a~point $x \in f(\Sigma)$ to another point of $f(\Sigma)$ is completely determined by how it acts on the ambient tangent space at~$x$~-- allowing us to identify ${\mathcal G}(f)$ with the groupoid $\mathrm{SO}[T_\Sigma M]$ of orientation-preserving orthonormal relative frames of $T_\Sigma M$. By a {\em relative frame} of a~vector bundle we mean an isomorphism between two fibres over possibly different basepoints. As $\Sigma $ is oriented, it has a well-defined unit normal, allowing us to identify $T_\Sigma M$ with the `thickened' bundle $T^+\Sigma :=T \Sigma \oplus ({\mathbb R} \times \Sigma)$, equipped with the obvious extension of the metric on $T \Sigma $ to an inner product. So ${\mathcal G}(f) \cong \mathrm{SO}[T^+\Sigma]$. Under this identification we see that the isotropy ${\mathcal G}(f)^2$ of $T\Sigma \subset T_\Sigma M$ is identified with $\mathrm{SO}[T\Sigma] \subset \mathrm{SO}[T^+\Sigma]$, the fundamental Lie groupoid associated with $\Sigma $ as a Riemannian manifold in its own right. Infinitesimalizing, we can identify the first two Lie algebroids in the invariant filtration $A(f) \supset A(f)^2 \supset A(f)^3\supset \cdots $ with the Lie algebroids $\mathfrak{so}[T^+ \Sigma] \supset \mathfrak{so}[T \Sigma]$. The same conclusion is drawn in Section \ref{bonnetsec2} by purely infinitesimal arguments.

Infinitesimal computations are straightforward, once one has the right model of the action algebroid ${\mathfrak g} \times M$ and a~corresponding description of its flat connection $\nabla $. As the action of $G \times M$ on $TM$ is faithful and preserves the metric, we have $G \times M \cong \mathrm{SO}[TM]$ and consequently
\begin{gather}
 {\mathfrak g} \times M \cong \mathfrak{so}[TM]
 \cong T M \oplus \mathfrak{so}(TM).\label{fpr}
\end{gather}
Here $\mathfrak{so}(TM)$ is the $\mathfrak{so}(n)$-bundle of skew-symmetric endomorphisms of the tangent bundle~$TM$ (the kernel of the anchor of $\mathfrak{so}[TM]$) and the second isomorphism is obtained with the help of the Levi-Cevita connection. In this model the connection $\nabla $ coincides with the canonical Cartan connection on~$\mathfrak{so}[TM]$, in the sense of~\cite{Blaom_06}, constructed for arbitrary Riemannian manifolds in~\cite{Blaom_12}.

Readers familiar with tractor bundles will recognise the model on the right of~\eqref{fpr} as the adjoint tractor bundle associated with $M$~\cite{Cap_Gover_02}. Those familiar with our work on infinitesimal geometric structures will recognise a model of the isotropy subalgebroid of the metric on~$M$ associated with a canonical Lie algebroid representation of the jet bundle $J^1 (TM)$ on $S^2\,T^*\!M$~\cite{Blaom_06,Blaom_12}. Similar models are known for other Klein geometries, and in particular for parabolic geometries~\cite{Cap_Slovak_09}, although the Lie algebroid structure of these models has been mostly ignored (notable exceptions are \cite{Armstrong_07,XuXiaomeng_14}).

The definitions of $\mathfrak{so}[T^+\Sigma] $ and $\mathfrak{so}[T \Sigma]$ depend on (and encode, up to scale) the metric that~$\Sigma $ inherits from $M$ but on no other aspect of the immersion $f$. The second fundamental form $\ii$ of the immersion enters into the present picture in two ways. Firstly, it enters indirectly in the form of $A(f)^3$, which we show is the isotropy subalgebroid of $\ii$, in the canonical representation of the Lie algebroid $A(f)^2 \cong \mathfrak{so}[T \Sigma] $ on the vector bundle $S^2\,T^*\Sigma$, of which $\ii$ is a section. The well-known fact, that $\Sigma \subset M $ has maximal symmetry when $\ii$ is a constant multiple of the inherited metric, follows immediately from this observation and a general characterisation of symmetries given in the next subsection. More generally, $A(f)^3$ is an {\em intransitive} Lie algebroid (and consequently does not appear in any formulation based on principal bundles, such as $G$-structures).

Secondly, we may recover $\ii$ directly as follows. The radical $\rad {\mathfrak g} $ of ${\mathfrak g} $ consists of all constant vector fields on $M ={\mathbb R}^n $, so that $\rad {\mathfrak g} \cong {\mathbb R}^n $. This identification transfers the standard inner product on ${\mathbb R}^n $ to an inner product $Q$ on $\rad {\mathfrak g}$ that is in fact $\Adjoint$-invariant, and we obtain an inner product on the trivial vector bundle $\rad E=\rad {\mathfrak g} \times \Sigma $ that we also denote by $Q$. It is not hard to see that $A(f)$, which we may identify with a subbundle of $E = {\mathfrak g} \times \Sigma $, intersects $\rad E$ in a subbundle of corank one. The {\em normal} for the immersion $f $ is the section $\xi$ of $\rad E$ defined, up to sign, by requiring $\xi $ to have $Q$-length one and be $Q$-orthogonal to $A(f) \cap \rad E$. The ambiguity in sign is resolved by requiring that the image ${\mathbf n}$ of $\xi $ under the natural projection $(\xi,x) \mapsto \xi(x) \colon E \rightarrow T_\Sigma M$~-- which is a unit length vector field orthogonal to $T \Sigma $ with respect to the metric on~$M$~-- have a direction consistent with the orientations on $\Sigma$ and $M $. Explicitly, we have at each point $x \in \Sigma$, $\xi (x)=(\xi'(x),x)$, where $\xi'(x)\in \rad {\mathfrak g} $ is the constant vector field on $M$ with $\xi'(x)(x)={\mathbf n}(x)$. Furthermore, $\nabla_u \xi \in A(f)$ for all $u \in T \Sigma $, where $\nabla $ is the canonical flat connection on $E={\mathfrak g} \times \Sigma $, and if fact
\begin{gather*}
 \ii(u,v) = -\langle\!\langle \# \nabla_u \xi,v\rangle\!\rangle,
\end{gather*}
where $\#$ denotes anchor and $\langle\!\langle \,\cdot\,,\,\cdot\,\rangle\!\rangle$ is the
metric on $\Sigma $.

Conversely, with ${\mathfrak g} $ continuing to denote the Lie algebra of Killing fields on ${\mathbb R}^n $ as above, suppose we are given an infinitesimal immersion $\omega \colon A \rightarrow {\mathfrak g} $, for some Lie algebroid $A$ over an oriented $(n-1)$-dimensional manifold $\Sigma$. Since $\omega $ is injective on fibres, it is an elementary observation that $T \Sigma $ can be regarded as a subbundle of $E/{\mathfrak h} $, where $E={\mathfrak g} \times \Sigma $ and ${\mathfrak h} $ is the kernel of the anchor of $A$. We define an inner product $Q$ on $\rad E$ as before and can show that the natural projection $E \rightarrow E/{\mathfrak h} $ restricts to an isomorphism $\rad E \rightarrow E/{\mathfrak h} $, pushing the inner product $Q$ to one on $E/{\mathfrak h} $. The restriction of this form to $T \Sigma \subset E/{\mathfrak h} $ is the {\df first fundamental form} of $\omega $, denoted $\langle\!\langle \,\cdot\,,\,\cdot\,,\rangle\!\rangle_\omega $, coinciding with the inherited metric on $\Sigma $ when $\omega $ is the logarithmic derivative of an immersion. After showing that $A$ must intersect $\rad E$ in a~corank-one subbundle, and explaining how $E/{\mathfrak h} $ can be naturally oriented, we are able to define the {\df normal} $\xi \in \Gamma(\rad E)$ of $\omega $ just as we did for immersions. The {\df second fundamental form} of~$\omega $, denoted $\ii_\omega $, is defined by
\begin{gather*}
 \ii_\omega(u,v) = -\langle\!\langle \# \nabla_u \xi,v\rangle\!\rangle_\omega,
\end{gather*}
and this coincides with the usual second fundamental form when $\omega $ is the logarithmic derivative of an immersion. Applying Theorem~\ref{infirmT}, we will obtain:
\begin{Theorem}[the abstract Bonnet theorem; proven in Section \ref{bonnetsec}]\label{introduceT} Assuming $\Sigma $ is simply-connected and $\rank A = (n+2)(n-1)/2$, any infinitesimal immersion $\omega \colon A \rightarrow {\mathfrak g} $ of $\Sigma $ into $M={\mathbb R}^n$ has, as primitive, an immersion $f \colon \Sigma \rightarrow {\mathbb R}^n $ with first and second fundamental forms $\langle\!\langle \,\cdot\,,\,\cdot\,,\rangle\!\rangle_\omega $ and $\ii_\omega $. The immersion $f $ is unique up to orientation-preserving isometries of~${\mathbb R}^n $.
\end{Theorem}

The classical Bonnet theorem for isometric immersions is obtained as a~corollary. In principle, Theorem~\ref{infirmT} could also be used to study monodromy obstructions to realising hypersurfaces in non-simply-connected cases but this is not attempted here.

A simpler illustration of the general theory, to planar curves, appears in Section \ref{invint}. An application to curves in the equi-affine plane is offered in Section~\ref{appendixA}.

\subsection*{Symmetries of a subgeometry}
Under our assumptions of constant rank, an invariant filtration $A(f)\supset A(f)^2\supset A(f)^3\supset \cdots $ must terminate at some Lie algebroid $A(f)^k$, $k\ge 1$. By construction, the flat connection $\nabla $ drops to a connection on $A(f)^k$ that is actually a Cartan connection (in the sense of~\cite{Blaom_06}). In particular, the space ${\mathfrak s}$ of $\nabla $-parallel sections of $A(f)^k$ is a Lie subalgebra of all sections, and acts infinitesimally on the submanifold $\Sigma $. We may view ${\mathfrak s}$ as a~subalgebra of $ {\mathfrak g} $ and ${\mathfrak s} $ has the dimension of the fibres of $A(f)^k$ because $\nabla $ is flat. Our Theorem~\ref{sssT} asserts that the infinitesimal action integrates to a~pseudogroup of transformations of $\Sigma $ consisting of restrictions to open subsets of $f(\Sigma) \subset M$ of symmetries of $M$.

\subsection*{Concluding remarks}
The results of \cite{Blaom_F} constitute a natural generalization of an elegant and widely applied result of \'Elie Cartan, which here delivers a~conceptual simplification in one of its main applications, the study
of subgeometry. Evidently there is an additional overhead, in the form of abstraction, in recovering known results in those relatively simple cases considered here. We expect the additional abstraction will be justified in more sophisticated applications.

A well-regarded and practical alternative to Cartan's moving frames which must be mentioned here, with diverse applications in mathematics and elsewhere, is due to Peter Olver and Mark Fels~\cite{Olver_Fels_98,Olver_Fels_99}; see also the survey~\cite{Olver_05}.

There are far simpler, coordinate-free proofs of the global Bonnet theorem for Euclidean geometry; see, e.g., \cite[Theorem~1.1]{Burstall_Calderbank_04}. However, the simpler methods do not generalize to arbitrary Klein geometries. In the case of the class of parabolic geometries there is a tractor bundle approach to subgeometry, involving a comparable degree of abstraction, due to Burstall and Calderbank \cite{Burstall_Calderbank_04}. In particular, this approach has been successfully applied to conformal geometry~\cite{Burstall_Calderbank_10}.

Possible extensions of the analysis initiated here include:

{\it Illustrations to other concrete Klein geometries.} Preliminary investigations suggest that hypersurfaces in affine geometry and conformal geometry are quite tractable, but it would be nice to obtain Bonnet-type theorems in novel cases also.

{\it Detailed descriptions of monodromy obstructions to hypersurface reconstruction in Riemannian and other geometries.} As we explain in Remark \ref{croonR}, there are generally {two} monodromy-type obstructions.

{\it Illustrations or novel applications to submanifolds of other codimension.} Even codimension zero (local diffeomorphisms) could be interesting to revisit, in particular with regard to monodromy obstructions.

{\it Descriptions of globalisations ${\mathcal G}(f)^j$ of the Lie algebroids $A(f)^j$, for $j\ge 3$.} It is not difficult to guess likely candidates for these (considering jets of immersions). Nice geometric interpretations of our infinitesimal analyses would likely follow. For example, why, from the present viewpoint, is the curvature of a curve in the Euclidean plane the radius of the circle of second-order tangency to the curve (Remark~\ref{planarR})?

{\it Unified description of curves $($or loops$)$ in planar geometries.} The examples investigated so far by the author, including the Euclidean and affine cases described here, suggest that in generic cases the last transitive Lie algebroid $A(f)^k$ in the invariant filtration is just a copy of the tangent bundle of the interval (or circle), delivering a right-inverse $s$ for the anchor of $A(f)\supset A(f)^k$. The `intrinsic' geometry of the curve is then (generically) a parallelism, equipping the curve with a canonical reparameterization $\tilde f$, up to scale; each polynomial invariant $Q$ on ${\mathfrak g} $ defines an ordinary function
$Q(\delta \tilde f(s(\partial /\partial t)))$ on the curve. Presumably these invariants suffice to describe the curve up to symmetry.

{\it Generalizations to submanifolds of curved geometries.} Suppose $M$ is a Cartan geometry, i.e., a Klein geometry deformed by curvature \cite{Sharpe_97}, and $f \colon \Sigma \rightarrow M$ an immersed submanifold. Then there exists a canonical Cartan connection $\nabla $ on an Atiyah Lie algebroid $A_M$ associated with the geometry \cite{Blaom_06,Blaom_16b}, a connection that generalizes the canonical flat connection on the action algebroids considered here. (A~more direct description of the pair $(A_M, \nabla)$ might also be possible, as we recall in the case of Riemannian geometry in Section~\ref{bonnetsec2}; see also~\cite{Blaom_12}.) One can use~$f$ to pull back $A_M$, in the category of Lie algebroids, to a subbundle $A(f)$ of $E$, where $E$ is the pullback of $A_M$ in the category of {\df vector bundles}. The connection $\nabla $ pulls back to a connection on $E$ and we may mimic the construction of the invariant filtration $A(f)\supset A(f)^2 \supset A(f)^3 \supset \cdots $ given above. This is still a filtration of Lie algebroids, and the fundamental representation of~$A(f)$ on~$E$ persists, as does the isotropy interpretation of transitive members of the filtration. For example, one can show that a~knowledge of this filtration in the case of a hypersurface of a~Riemannian manifold suffices to describe its inherited metric and second fundamental form, up to a choice of scale.

{\it Parabolic geometry and Courant algebroids.} As one referee suggests, in the special case of parabolic geometries, it may be appropriate to consider the underlying Courant algebroid structure \cite{Armstrong_07,XuXiaomeng_14}, rather than the Lie algebroid one. We note, however, that no analogue of Lie~II for Courant algebroids (and hence of Theorem~\ref{infirmT}) is currently known.

\subsection*{Bracket convention and notation} Throughout this article, brackets on Lie algebras and Lie algebroids are defined using {\em right}-invariant vector fields. We reserve the symbol $G$ for Lie groups and ${\mathfrak g} $ for Lie algebras. The symbols ${\mathcal G} $, ${\mathcal G}(f)$, and ${\mathcal G}(f)^2$ denote Lie groupoids; $A$, $A^2$, $A(f)$, $A(f)^2$, etc., denote Lie algebroids. The kernel of the anchor $\#$ of $A$ or $A(f)$ is a Lie algebra bundle denoted ${\mathfrak h} $, that of $A^2$ or $A^2(f)$ is denoted by ${\mathfrak h}^2$, etc. We reserve $E$ for the trivial ${\mathfrak g} $-bundle ${\mathfrak g} \times \Sigma $ over $\Sigma$.

The frame groupoid of a vector bundle $B $ will be denoted $\mathrm{GL}[B ]$ and the corresponding subgroupoid of automorphisms of fibres of $B $ by $\mathrm{GL}(B )$. So, an element of $\mathrm{GL}[B ]$ is an isomorphism (relative frame) $B |_{m_1} \rightarrow B |_{m_2}$, becoming an element of $\mathrm{GL}(B )$ if $m_1=m_2$. Similar notation applies to the subgroupoids and Lie subalgebroids of frame groupoids.

\section{Symmetries and infinitesimal immersions}\label{infch}
In formulating Theorem \ref{infirmT}, several notions were introduced on which we now elaborate.

\subsection*{Simple Klein geometries}
Let $M$ be a Klein geometry with transitively acting Lie group
$G$. Let $G_{m_0}^\circ $ denote the connected component of the
isotropy $G_{m_0}$ at $m_0 \in M$. As $G_{m_0}^\circ$ is {\em
 path}-connected, $N_G(G_{m_0})\subset N_G(G_{m_0}^\circ)$. Here
$N_G(H)$ is the normaliser of $H \subset G$ in $G$.
\begin{Definition}\label{weakD}
 We say the isotropy groups of the $G$ action are {\df weakly
 connected} if for some (and hence any) $m_0 \in M$, we have
 $N_G(G_{m_0})=N_G(G_{m_0}^\circ)$.
\end{Definition}
\noindent%
As already mentioned, $M$ is called a {\df simple Klein geometry} in
this case.

\subsection*{Symmetries of a Klein geometry}
For reasons explained in \cite{Blaom_F}, we allow symmetries of Klein geometries to include more than the left translations $m \mapsto g \cdot m$, $g \in G$:
\begin{Definition}\label{symmetriesD}
 Let $M$ be a simple Klein geometry and $G$ the transitively acting Lie group. Then a {\df symmetry} of $M$ is any diffeomorphism $\phi \colon M \rightarrow M$ for which there exists some $l \in G$ such that $\phi(g \cdot m) = lgl^{-1} \cdot \phi(m)$ for all $g \in G,\, m \in M$.
\end{Definition}

The symmetries of $M$ form a Lie group, henceforth denoted $\automorphism(M)$. If we take $M=G$ (with~$G$ acting on itself from the left) then $\automorphism(M)$ consists of all left {\em and right} translations. However, $\automorphism(M)$ is frequently not much larger than the group of left translations, in applications of interest to geometers. For example, to the group~$G$ of orientation-preserving isometries of a Riemannian space form
${\mathbb R}^n$, $S^n$ or ${\mathbb H}^n $, one has only to add the orientation-{\em reversing} isometries to obtain the full symmetry group $\automorphism(M)$, and then only in the case of even-dimensional spheres. If, in these examples, we instead take $G$ to be the full group of isometries (in which case $M$ is still a~simple Klein geometry) then we even have $\automorphism(M)=G$. For further examples, see \cite[Examples~2.3]{Blaom_F}.

\subsection*{Infinitesimal immersions and their synthesis}
Let $M$ be a Klein geometry, $G$ the transitively acting Lie group, and $\Sigma $ a fixed connected smooth manifold. Then a Lie algebroid morphism $\omega \colon A \rightarrow {\mathfrak g} $ is an {\em infinitesimal immersion} of $\Sigma $ into $M$ if $A $ has base $\Sigma $ and $\omega$ satisfies the following properties, satisfied by the logarithmic derivatives of bona fide immersions $f \colon \Sigma \rightarrow M$:
\begin{enumerate}\itemsep=0pt
\item[\bf I1.] $A$ is transitive.
\item[\bf I2.] $\omega $ is injective on fibres.
\item[\bf I3.] For some point $x_0 \in \Sigma $ (and hence {\em any} point $x_0 \in \Sigma $; see \cite[Theorem 2.4]{Blaom_F}) there exists $m_0 \in M$ such that $x_0\xrightarrow{\omega}m_0$.
\item[\bf I4.] $\dim {\mathfrak g}-\rank A = \dimension M - \dimension \Sigma$.
\end{enumerate}
Here $x_0\xrightarrow{\omega}m_0$ is shorthand for the condition $\omega(A_{x_0}) \subset {\mathfrak g}_{m_0}$, where $A_{x_0}$ denotes the isotropy algebra of $A$ at ${x_0} \in \Sigma $, i.e., the kernel of the restriction of the anchor to the fibre~$A|_{x_0}$, while~${\mathfrak g}_{m_0}$ denotes the usual isotropy (stabiliser) at~$m_0$ of the infinitesimal action of ${\mathfrak g} $ on $M$. Axiom~I4 ensures that $A$ has the largest rank consistent with Axioms~I1 and~I2, in which case Axiom~I3 actually implies
\begin{gather}
 \omega(A_{x_0})= {\mathfrak g}_{m_0}. \label{worm}
\end{gather}

In the language of \cite{Blaom_F}, $\omega $ is an infinitesimal immersion if it is a maximal generalized Maurer--Cartan form that is injective on fibres.\footnote{$\omega $ is a {\df generalized Maurer--Cartan form} if it satisfies Axioms~I1, I3 and a weakened form of Axiom~I2: the restriction of $\omega $ to every isotropy algebra $A_x$ should be injective, $x \in \Sigma $.}

Now suppose $\omega \colon A \rightarrow {\mathfrak g} $ is an infinitesimal immersion of $\Sigma $ into $M$. By virtue of Axiom~I2, we may use $\omega $ to identify $A$ with a subbundle of the trivial bundle $E := {\mathfrak g} \times \Sigma \rightarrow \Sigma $. This bundle comes equipped with a canonical flat connection $\nabla $ and a~canonical $\nabla$-parallel Lie bracket $\{\,\cdot \,,\,\cdot\,\}$ on its space of sections making it into a ${\mathfrak g} $-bundle: $\{X,Y\}(x):=[X(x),Y(x)]_{\mathfrak g}$. Applying the characterisation of Lie algebroid morphisms in terms of connections (see, e.g., \cite[Proposition~4.1.9, p.~154]{Mackenzie_05}) we deduce the {\df $($generalized$)$ Maurer--Cartan equation},
\begin{gather}
 [X,Y] = \nabla_{\#X}Y - \nabla_{\#Y}X +\{X,Y\}, \label{mce}
\end{gather}
holding for all sections $X$, $Y$ of $A$. Implicit in this formula is that the right-hand side is a section of $A \subset E$.

The preceding observations have a readily established converse. By a {\df ${\mathfrak g}$-bundle} let us mean a~Lie algebra bundle modelled on ${\mathfrak g} $, in the sense of, e.g.,~\cite{Mackenzie_05}. Then:
\begin{Proposition}\label{croonP} Let $E$ be a ${\mathfrak g}$-bundle over a simply-connected manifold $\Sigma $, equipped with a flat connection~$\nabla $, with respect to which the bracket on $E$ is parallel. Then:
\begin{enumerate}\itemsep=0pt
 \item[$1.$] 
 The space ${\mathfrak g}'$ of $\nabla $-parallel sections of $E$ is a Lie subalgebra of $\Gamma(E)$ isomorphic to ${\mathfrak g} $ and the map $(\xi,x)\mapsto \xi(m) \colon {\mathfrak g}' \times \Sigma \rightarrow E$ is a vector bundle isomorphism.
 \item[$2.$] 
 For any Lie algebroid $A$ over $\Sigma $, realised as subbundle of $E$ in such a way that the Maurer--Cartan equations~\eqref{mce} hold, the composite $A \hookrightarrow E \rightarrow {\mathfrak g}' \cong {\mathfrak g}$ is an injective Lie algebroid morphism $\omega \colon A \rightarrow {\mathfrak g} $. Here $E \rightarrow {\mathfrak g}'$ is the projection determined by the isomorphism $E\cong {\mathfrak g}' \times \Sigma$ in~{\rm (1)}.
\end{enumerate}
\end{Proposition}
Equipped with the preceding proposition, we posit the following general principle, basic to later illustrations of the general theory:
\begin{constructionPrinciple}\label{constructionPrinciple} To construct an infinitesimal immersion $\omega \colon A\rightarrow {\mathfrak g} $, attempt to realize the Lie algebroid $A$ as a subbundle of a $\mathfrak g$-bundle $E$ equipped with a flat connection $\nabla $, with respect to which the bracket on $E$ is parallel, and chosen such that the Maurer--Cartan equations \eqref{mce} hold. Candidates for $A$, $E$ and $\nabla $ are suggested by computing the logarithmic derivatives of smooth immersions.
\end{constructionPrinciple}
\begin{Remark}\label{croonR} If we drop the proposition's simple connectivity hypothesis then there is a global obstruction to synthesising a Lie algebroid morphism $\omega \colon A \rightarrow {\mathfrak g} $ in this way, namely the mono\-dromy of $\nabla $. So, in general, there are {\em two} monodromy obstructions to reconstructing immersions from infinitesimal data: the monodromy of a connection~-- necessary for constructing an infinitesimal immersion $\omega $~-- and the monodromy of $\omega $ itself (see~\cite{Blaom_F}).
\end{Remark}

For a first illustration of the Construction Principle, see Section~\ref{invint} under `Planar curves'.

\section{Invariants}\label{invint}

\subsection*{Morphisms between infinitesimal immersions}
Fixing a Klein geometry $M$ with transitively acting group $G$, and a connected manifold $\Sigma $, we will collect all infinitesimal immersions of $\Sigma $ into $M$ into the objects of a category $\mathbf{Inf}$. In this category a morphism $\omega_1 \rightarrow \omega_2$ between objects $\omega_1 \colon A_1 \rightarrow {\mathfrak g} $ and $\omega_2 \colon A_2 \rightarrow {\mathfrak g} $ consists of a Lie algebroid morphism $\lambda \colon A_1 \rightarrow A_2$ covering the identity on $\Sigma $ and an element $l \in G$ such that the following diagram commutes:
\begin{gather*}
 \begin{CD}
 A_1 @>{\omega_1}>> {\mathfrak g} \\
 @V{\lambda}VV @VV{\Adjoint_l}V \\
 A_2 @>{\omega_2}>> {\mathfrak g}.
 \end{CD}
\end{gather*}
Axioms I2 and I4 ensure that all morphisms in $\mathbf{Inf}$ are, in fact, isomorphisms. In this language, a smooth map $f \colon \Sigma \rightarrow M$ is a primitive of an infinitesimal immersion $\omega $ if and only if $\omega$ and $\delta f $ are isomorphic in $\mathbf{Inf}$. These abstractions are further justified by the following:

\begin{Theorem}[\cite{Blaom_F}]\label{morphismsT}
 Suppose $M$ is a simple Klein geometry and let $f_1, f_2 \colon \Sigma \rightarrow M$ be smooth maps. Then there exists an isomorphism $\delta f_1 \cong \delta f_2$ in $\mathbf{Inf}$ if and only if there exists $\phi \in \automorphism(M)$ such that $f_2 = \phi\, \circ f_1$.
\end{Theorem}

In particular, smooth immersions $f_1,f_2 \colon \Sigma \rightarrow M$ agreeing up to a symmetry of $M$ have isomorphic logarithmic derivatives.

\subsection*{Invariants defined}
While the magnanimous reader will generally understand the term `invariant' without further comment, to mitigate possible confusion we state here formal definitions sufficient for the sequel. Stronger definitions are possible, but only at the cost of further abstraction we prefer to avoid.

Suppose ${\mathcal C} $ is a collection of smooth immersions $f \colon \Sigma \rightarrow M$ closed under the action of $\automorphism(M)$. Then a map $f \mapsto Q_f$ from ${\mathcal C} $ to some category will be called an {\df invariant} for ${\mathcal C} $ if, for all $\phi \in \automorphism(M)$, there exists an isomorphism $Q_{\phi \circ f} \cong Q_f$. An {\df invariant} for subcategory ${\mathcal D} $ of $\mathbf{Inf}$ will be any functor
$\omega \mapsto Q_\omega $ from ${\mathcal D} $ into another category. A trivial corollary of the preceding theorem is:
\begin{Proposition}\label{rhinoP}
Every invariant for infinitesimal immersions delivers an invariant for smooth immersions. Specifically, with ${\mathcal C} $ and ${\mathcal D} $ as above, and assuming $\delta ({\mathcal C}) \subset {\mathcal D} $, every invariant $\omega \mapsto Q_\omega $ for ${\mathcal D}\subset \mathbf{Inf} $ delivers an invariant $f \mapsto Q_{\delta f}$ for the collection of smooth immersions ${\mathcal C} $.
\end{Proposition}

\subsection*{Invariants defined by polynomials}
Let $M$ be a Klein geometry with transitively acting group $G$, and $Q \colon {\mathfrak g} \rightarrow {\mathbb K} $ a polynomial invariant under the adjoint representation of $G$; here ${\mathbb K}$ is ${\mathbb R} $ or ${\mathbb C} $. For each infinitesimal immersion $\omega \colon A \rightarrow {\mathfrak g} $ we define
\[Q_\omega =Q \circ \omega \colon \ A(f)\rightarrow {\mathbb K}\]
and immediately obtain:
\begin{Proposition}\label{polynomialP}
 Let $\omega \colon A \rightarrow {\mathfrak g} $ and $\omega' \colon A' \rightarrow {\mathfrak g} $ be isomorphic infinitesimal immersions and $\lambda \colon A \rightarrow A'$ a corresponding isomorphism. Then the following diagram commutes
 \[
 \begin{tikzcd}
 A \arrow{d}{\lambda } \arrow{r}{Q_{\omega}} & {\mathbb K} \\
 A'.\arrow[swap]{ru}{Q_{\omega'}} &
 \end{tikzcd}
 \]
 In particular, $\omega \mapsto Q_\omega$ is an invariant for $\mathbf{Inf}$.
\end{Proposition}
Similarly, for any $\Adjoint$-invariant subalgebra ${\mathfrak g}' \subset {\mathfrak g} $ (e.g., the commutator $[{\mathfrak g},{\mathfrak g}]$ or radical of ${\mathfrak g} $) the pre-image $A':=\omega^{-1}({\mathfrak g} ')$ is Lie subalgebroid of $A$, assuming it has constant rank (because~$\omega $ is a Lie algebroid morphism) and is an invariant for $\mathbf{Inf}$. Any $\Adjoint$-invariant polynomial on ${\mathfrak g}'$ furnishes yet
another invariant.

\begin{Remark}\label{polynomialR}
In applying the Construction Principle \ref{constructionPrinciple} it is important to distinguish between polynomials on ${\mathfrak g} $ that are merely invariant under the adjoint action ({\df inner} invariants) and those invariant under arbitrary Lie algebra automorphisms ({\df outer} invariants) for the following reason. In Proposition~\ref{croonP}(2) the infinitesimal immersion $\omega \colon A \rightarrow {\mathfrak g} $ is only determined up to an automorphism of ${\mathfrak g} $, because of the choice of isomorphism ${\mathfrak g}' \cong {\mathfrak g} $. Therefore, given an $\Adjoint$-invariant polynomial $Q$ on ${\mathfrak g} $, the corresponding map $Q_\omega \colon A \rightarrow {\mathbb K} $ depends on the choice made, unless $Q$ is an outer invariant (e.g., if $Q$ is the Killing form). Concrete illustrations given later will make this point clearer.
\end{Remark}

Note that an $\Adjoint$-invariant polynomial $Q \colon {\mathfrak g} \rightarrow {\mathbb K} $ defines, in a trivial way, a map $Q \colon E \rightarrow {\mathbb K} $ on the trivial bundle $E={\mathfrak g} \times \Sigma $, which will be called an {\df invariant polynomial on $E$} and denoted~$Q$ also. Then, when we use an infinitesimal immersion $\omega \colon A \rightarrow {\mathfrak g}$ to identify $A$ with a subbundle of $E={\mathfrak g} \times \Sigma $, $Q_\omega $ is the composite of the inclusion $A \hookrightarrow E$ with $Q \colon E \rightarrow {\mathbb K} $.

\subsection*{The fundamental representation}
Let $\omega \colon A \rightarrow {\mathfrak g} $ be an infinitesimal immersion of $\Sigma$ into $M$ and identify the Lie algebroid $A$ with a subbundle of $E={\mathfrak g} \times \Sigma $. Continuing to let $\nabla $ denote the canonical flat connection on $E$, one has the following representation $\bar\nabla $ of the Lie algebroid $A$ on~$E$:
\begin{gather}
 \bar\nabla_X \xi = \nabla_{\#X} \xi + \{X, \xi\}.\label{hsq}
\end{gather}
That $\bar\nabla $ is indeed a representation follows from the Maurer--Cartan equations \eqref{mce}, the fact that $\{\,\cdot\,,\,\cdot\,\}$ is $\nabla $-parallel, and the Jacobi identity for $\{\,\cdot\,,\,\cdot\,\}$. We shall call this representation the {\df fundamental representation}. We note that the Maurer--Cartan equations also allow us to write
\begin{gather}
 \bar\nabla_XY=\nabla_{\#Y} X + [X,Y],\label{aall}
\end{gather}
in the special case that $X$ and $Y$ are both sections of $A$.
\begin{Lemma}\label{fffL}
 Any $\Adjoint$-invariant polynomial $Q$ on $E$ is invariant with respect to the fundamental representation. The radical $\rad E=\rad{\mathfrak g} \times \Sigma $, the commutator $\{E,E\}$, and elements of the upper or lower derived series of $E $, are all invariant under the fundamental representation.
\end{Lemma}
\begin{proof}
For the first statement, suppose $Q$ is quadratic, the other cases being similar. View $Q$ as a symmetric bilinear form on $E$. Then, because $\Adjoint$-invariance implies $\adjoint$-invariance, we have $Q(\{X, \xi_1\}, \xi_2) + Q(\xi_1, \{X,\xi_2\})$ for any $X, \xi_1, \xi_2 \in \Gamma(E)$. It is trivial that $Q$ is $\nabla $-parallel, and so not hard to see from~\eqref{hsq} that $\bar\nabla_XQ=0 $ for any $X \in \Gamma(A)$. The second statement follows from the fact that the bracket $\{\,\cdot\,,\,\cdot\,\}$ is invariant with respect to the fundamental representation, which follows from its $\nabla $-invariance and the Jacobi identity.
\end{proof}

The fundamental representation gives an alternative characterization of the invariant filtration to be described next.

\subsection*{The invariant filtration}
In general $A \subset E $ is not $\nabla $-invariant (or invariant under the fundamental representation $\bar\nabla$). We are therefore led to define, for each $x \in \Sigma $,
\begin{gather*}
 A^2|_x =\{X(x)\suchthat\text{$X \in \Gamma(A)$ and $\nabla_uX \in A$ for all $u \in T_x \Sigma$}\},
\end{gather*}
and, assuming the dimension of $A^2|_x$ is independent of $x$, obtain a vector subbundle $A^2 := \amalg_{x \in \Sigma}A^2|_x\subset A$. Tacitly assuming rank-constancy at each stage, we inductively define, for each $j \ge 3$, a subbundle $A^j \subset A^{j-1}$ with fibres given by
\begin{gather*}
 A^{j}|_x =\big\{X(x)\suchthat\text{$X \in \Gamma\big(A^{j-1}\big)$ and $\nabla_uX \in A^{j-1}$ for all $u \in T_x \Sigma$}\big\},
\end{gather*}
and obtain a filtration
\[ A \supset A^2 \supset A^3 \supset A^4 \supset \cdots \]
of vector subbundles, which are in fact subalgebroids, by the proposition below. This filtration, determined by the infinitesimal immersion $\omega $, will be called the associated {\df invariant filtration}, for one readily proves the following:
\begin{Lemma}\label{fff2L}
 Let $\omega$ and $\omega'$ be isomorphic Maurer--Cartan forms and $L \colon A \rightarrow A' $ an isomorphism, as in Proposition~{\rm \ref{polynomialP}}. Then $L(A^j) = (A')^j$ for all $j$.
\end{Lemma}

In particular, the restriction to $A(f)^j$ of a invariant polynomial $Q \colon E \rightarrow {\mathbb K} $ of a smooth map $f \colon \Sigma \rightarrow M$ is an invariant of $f$, generally `finer' than its restriction to~$A(f)$.

\begin{Proposition}\label{fff2P} For all $j \ge1$, $A^j \subset A$ is a subalgebroid. In the case that $A^j$ is transitive, $A^{j+1}$ is the isotropy subalgebroid of $A^j \subset E$ under the restriction of the fundamental representation to a representation of $A^j$.
\end{Proposition}
\begin{proof} For a proof of the first statement by an induction on $j$, suppose $A^j \subset A$ is a subalgebroid and let $X$ and $Y$ be sections of $A^{j+1}$, so that $\nabla_UX$ and $\nabla_UY$ are sections of $A^j$ for all vector fields $U$ on $\Sigma $. We must show $\nabla_U[X,Y]$ is a section of $A^j$ also. To this end, replace $X$ in the Maurer--Cartan equations \eqref{mce}, with $\nabla_UX$ to conclude
 \begin{gather}
 \nabla_{\#Y}\nabla_U X - \{\nabla_UX, Y\}\text{~is a section of $A^j$}\label{rone}.
 \end{gather}
 Similarly, replace $Y$ in the Maurer--Cartan equations with $\nabla_UY$ to conclude
 \begin{gather}
 \nabla_{\#X} \nabla_U Y+ \{X, \nabla_U Y\}\text{~is a section of $A^j$}\label{rtwo}.
 \end{gather}
Finally, apply $\nabla_U$ to both sides of the Maurer--Cartan equations and, appealing to the flatness of $\nabla $ and the $\nabla $-invariance of $\{\,\cdot \,,\,\cdot\,\}$, show that
 \begin{gather*}
 \nabla_U[X,Y]= (\nabla_{\#X} \nabla_U + \{X, \nabla_U Y\}) -(\nabla_{\#Y}\nabla_U X - \{\nabla_UX, Y\})\\
\hphantom{\nabla_U[X,Y]=}{} +\nabla_{[U,\#Y]}X - \nabla_{[U,\#X]}Y.
 \end{gather*}
 With the help of \eqref{rone} and \eqref{rtwo} we now see that $\nabla_U[X,Y]$ must be a section of $A^j$.

 The second statement in the proposition follows easily enough from~\eqref{aall}.
\end{proof}

\subsection*{Symmetries of infinitesimal immersions and subgeometries}
A connection $\nabla $ on an arbitrary Lie algebroid $A$ is a {\df Cartan connection} if it suitably respects the Lie algebroid structure of $A$. The canonical flat connection $\nabla $ on the action algebroid ${\mathfrak g} \times M$ is such a connection. Conversely, any Lie algebroid equipped with a {\em flat} Cartan connection is called by us a {\df twisted Lie algebra action}, for it is locally an action algebroid \cite{Blaom_06}, and globally an action algebroid with `monodromy twist'~\cite{Blaom_13}. In the same way that a Lie algebra action integrates to a local Lie group action, so, more generally, every twisted Lie algebra action integrates to a~{\em pseudoaction}, a geometric object encoding the pseudogroup of transformations generated by the flows of the (locally defined) infinitesimal generators~\cite{Blaom_16}.

Let $\{A^j\}_{j\ge 1}$ be the invariant filtration defined by an infinitesimal immersion $\omega \colon A \rightarrow {\mathfrak g} $. Notice that if $A^{k+1}= A^k$, for some $k \ge 1$, then $A^j = A^k$ for all $j\ge k$~-- i.e., the filtration has stabilised at $A^k$. In our finite-dimensional setting all invariant filtrations must stabilise eventually (under our constant rank assumption) although $A^k =0 $ is typical.
\begin{Theorem}\label{sssT} Suppose $A^{k+1}=A^k$. Then the canonical flat connection $\nabla $ on $E$ restricts to a flat Cartan connection on $A^k$, i.e., defines a twisted Lie algebra action on $\Sigma$ which is, in fact, an ordinary action by a subalgebra ${\mathfrak s}\subset {\mathfrak g} $. In the case that $\omega $ is a logarithmic derivative of a~smooth immersion $f \colon \Sigma \rightarrow M$, elements of the corresponding pseudogroup of transformations on~$\Sigma $ correspond to restrictions to open subsets of $f(\Sigma) \subset M$ of left translations $m \mapsto g \cdot m$ in~$M$, $g \in G$.
\end{Theorem}
\begin{proof} If $A^{k+1}=A^k$, then $\nabla $ restricts to a connection on $A^k$ by the definition of $A^{k+1}$. Now the restricted connection~-- also denoted~$\nabla $~-- is Cartan if its cocurvature vanishes. But, from the definition of cocurvature (see~\cite{Blaom_06}) follows the formula
 \begin{gather*}
 \cocurvature \nabla (X,Y)U=\nabla_U\torsion \bar\nabla\,(X,Y) + \curvature \nabla (U, \#X)Y -\curvature \nabla (U, \#Y) X,
 \end{gather*}
where $\torsion \bar\nabla (X,Y) = \nabla_{\#Y}X - \nabla_{\#X}Y + [X,Y]$.
In the present case $\curvature \nabla =0$, while $\torsion \bar\nabla(X,Y)$ $= \{X,Y\}$, on account of the Maurer--Cartan equations \eqref{mce}. So the cocurvature vanishes.

Since $\nabla $ is a flat Cartan connection and $\Sigma $ is simply-connected, the subspace ${\mathfrak s} \subset \Gamma(A^k)$ of $\nabla$-parallel sections is a Lie subalgebra acting infinitesimally on $\Sigma $ according to $\xi^\dagger(x)=\#(\xi(x))$ \cite[Theorem~A]{Blaom_06}. Trivially, every such section is also $\nabla $-parallel as a section of $E={\mathfrak g} \times \Sigma $, and so is nothing but the restriction to $\Sigma $ of an infinitesimal generator of the action of ${\mathfrak g}$ on $M$. In particular, its local flows are restrictions to open subsets of $\Sigma $ of left translations $m \mapsto g \cdot m$ in~$M$, $g \in G$.
\end{proof}

\subsection*{The invariant filtration in terms of isotropy}
The isotropy characterization of Proposition \ref{fff2P} of the Lie algebroids appearing in an invariant filtration can be described more concisely. While $A \subset E$ is not generally invariant with respect to the fundamental representation, the kernel ${\mathfrak h} \subset A$ of its anchor {\em is} invariant, by the identity~\eqref{aall}. It follows that the fundamental representation drops to a representation of $A $ on $E/{\mathfrak h} $. Under the canonical identification $A/{\mathfrak h} \cong T \Sigma $, which we have by virtue of the transitivity of $A$, we obtain a canonical inclusion of $T \Sigma $ into $ E/{\mathfrak h} $. Evidently, $A^2 \subset A$ must coincide with the isotropy of $T \Sigma \subset E/{\mathfrak h} $ with respect to this representation. One continues in a similar fashion for the rest of the filtration and obtains, with $A^0=E$ and $A^1=A$:

\begin{Proposition}\label{cP}
 Assume $A^j$ is transitive, $j\ge 1$, and let ${\mathfrak h}^j \subset A^j$ denote the kernel of its anchor. Then the fundamental representation restricts to a representation of $A^j$ that leaves $A^{j-1}\subset E$ and ${\mathfrak h}^j \subset A^{j-1} $ invariant, dropping to a representation of $A^j$ on $A^{j-1}/{\mathfrak h}^j$. Moreover, $A^{j+1}$ is the isotropy of $T \Sigma \subset A^{j-1}/{\mathfrak h}^j$ when one views $T \Sigma$ as a subbundle of $ A^{j-1}/{\mathfrak h}^j$ using the isomorphism $T \Sigma \cong A^j/{\mathfrak h}^j$ determined by the anchor of $A_j$. In particular, $A^{j+1}$ acts canonically on $T \Sigma $.
\end{Proposition}
\begin{Remark}\label{cR} In the special case that $\omega $ is the logarithmic derivative of an immersion $f \colon \Sigma \rightarrow M$, $E/{\mathfrak h} $ has a canonical identification with the pullback $T_\Sigma M$ of $TM$ because the morphism $E \rightarrow T_\Sigma M$ sending $(\xi, m)$ to $\xi^\dagger(m)$ is surjective. So the fundamental representation induces a~representation of~$A(f)$ on $T_\Sigma M$ and $A(f)^2$ is the isotropy of $T \Sigma \subset T_\Sigma M$ under this
 representation.
\end{Remark}

\subsection*{Orientability of the abstract model of $\boldsymbol{T_\Sigma M}$}
Of course if $M$ is orientable, then so is $T_\Sigma M$. In general, however, the orientability of $E/{\mathfrak h} $ is a~subtle question. The issue is settled in the simplest scenario as follows:
\begin{Proposition}\label{orientedP} Let $\omega \colon A \rightarrow {\mathfrak g} $ be an infinitesimal immersion and ${\mathfrak h} $ the kernel of the anchor $\# \colon A \rightarrow T \Sigma $. Assume that $M$ is oriented and that the action of~$G$ on $M$ is orientation-preserving. Suppose, moreover, that for some $m_0 \in M$, we have $N_G(G_{m_0})=G_{m_0}$ and that~$G_{m_0}$ is connected. Then the orientation on $M$ determines a natural orientation of $E/{\mathfrak h}$ and this orientation is an invariant of $\omega$. If $\omega $ is the logarithmic derivative of an infinitesimal immersion, then the orientation of~$E/{\mathfrak h}$ coincides with that of $T_\Sigma M$ under the canonical identification \smash{$E/{\mathfrak h} \cong T_\Sigma M$}.
\end{Proposition}
\begin{proof} Let $x \in \Sigma $. Then by Axiom I3 and its consequence \eqref{worm}, there exists $m \in M$ such that the fibre ${\mathfrak h}|_x$ coincides with ${\mathfrak g}_m \times \{x\} \subset E$, under our tacit identification of ${\mathfrak h} $ with a subbundle of $E={\mathfrak g} \times \Sigma $. So $(E/{\mathfrak h})|_x \cong {\mathfrak g}/{\mathfrak g}_m$. For some $g \in G$ we have ${\mathfrak g}_m = \Adjoint_g({\mathfrak g}_{m_0})$ (pick $g$ to satisfy $m=g \cdot m_0$). The morphism $\Adjoint_g \colon {\mathfrak g} \rightarrow {\mathfrak g} $ drops to an isomorphism ${\mathfrak g}/{\mathfrak g}_{m_0} \rightarrow {\mathfrak g}/{\mathfrak g}_m$, also denoted~$\Adjoint_g$ below. Now the action of $G$ determines an canonical isomorphism $T_{m_0}M \cong {\mathfrak g}/{\mathfrak g}_{m_0}$, so we obtain a sequence of isomorphisms
\begin{gather}
 T_{m_0}M \cong {\mathfrak g}/{\mathfrak g}_{m_0} \xrightarrow{\Adjoint_g} {\mathfrak g}/{\mathfrak g}_m\cong (E/{\mathfrak h})|_x.\label{code}
\end{gather}
We use this isomorphism to transfer the orientation of $M$ to an orientation of $(E/{\mathfrak h})|_x$. It remains to show that this orientation is independent of the choices of $m \in M$ and $g \in G$. Indeed, suppose that ${\mathfrak g}_m={\mathfrak g}_{m'}$ and $\Adjoint_{g'}({\mathfrak g}_{m_0})= {\mathfrak g}_{m'}$, for some $m' \in M$ and $g' \in G$. Then $\Adjoint_{g'}({\mathfrak g}_{m_0})={\mathfrak g}_m$ and so $\Adjoint_{k}({\mathfrak g}_{m_0})={\mathfrak g}_{m_0}$, where $k=g'g^{-1}$. Since $G_{m_0} $ is connected, this implies $kG_{m_0}k^{-1}=G_{m_0}$, i.e., $k \in N_G(G_{m_0})=G_{m_0}$. With our new choices the isomorphism in \eqref{code} is replaced by the composite
\begin{gather*}
 T_{m_0}M \cong {\mathfrak g}/{\mathfrak g}_{m_0}\xrightarrow{\Adjoint_k}{\mathfrak g}/{\mathfrak g}_{m_0} \xrightarrow{\Adjoint_g} {\mathfrak g}/{\mathfrak g}_m\cong (E/{\mathfrak h})|_x.
\end{gather*} Because $k \in G_{m_0}$, this composite coincides with the composite
\begin{gather}
 T_{m_0}M \xrightarrow{TL_k}T_{m_0}M \cong {\mathfrak g}/{\mathfrak g}_{m_0} \xrightarrow{\Adjoint_g} {\mathfrak g}/{\mathfrak g}_m\cong (E/{\mathfrak h})|_x, \label{hoppy}
\end{gather}
where \looseness=1 $TL_k$ is the tangent lift of $L_k(n):=k \cdot n$, $n \in M$. As the action of $G$ on $M$ is orientation-preserving, whether we declare an orientation on $(E/{\mathfrak h})|_x$ using the isomorphism \eqref{code} or \eqref{hoppy} therefore makes no difference. The claims made regarding this orientation follow by construc\-tion.
\end{proof}

\subsection*{The globalisation of $\boldsymbol{A(f)^2}$}
We now inject some observations which are not essential to the computation of invariant filtrations, but which help to anticipate their outcomes in the case of a logarithmic derivative. Let $f \colon \Sigma \rightarrow M$ be an injective immersion. By definition, we may view $A(f)$ as the Lie algebroid of the pullback ${\mathcal G}(f) $ by $f$ of the action groupoid $G \times M$, given by~\eqref{ssaa}. The tangent-lifted action of $G$ on~$TM$ determines an action of the Lie groupoid $G \times M$ on $TM$, which pulls back to an action of ${\mathcal G}(f)$ on $T_\Sigma M$.
\begin{Proposition}\label{nungP} Let ${\mathcal G}(f)^2 \subset {\mathcal G}(f)$ denote the isotropy of $T \Sigma \subset T_\Sigma M$ under the representation of ${\mathcal G}(f)$ on $T_\Sigma M$. Then $A(f)^2 \subset A(f)$ is the Lie algebroid of ${\mathcal G}(f)^2$.
\end{Proposition}
\begin{proof} By Remark \ref{cR}, it suffices to show that the representation of $A(f)$ on $T_\Sigma M$ determined by the fundamental representation is nothing more than the infinitesimalization of the representation of ${\mathcal G}(f)$ on $T_\Sigma M$ described above. To this end, note first that the tangent action of $G \times M$ on $TM$ infinitesimalizes to a representation $\bar\nabla $ of ${\mathfrak g} \times M$ on $TM$ given by
 \begin{gather}
 \bar\nabla_X U = \# \nabla_U X + [\#X, U],\label{stardust}
 \end{gather}
 where $\nabla $ is the canonical flat connection on ${\mathfrak g} \times M$ and $\#$ the anchor map, $\# (\xi,m)=\xi^\dagger(m)$. If $U = \# Y$ for some section $Y$ of ${\mathfrak g} \times M$ then we may rewrite~\eqref{stardust} as
 \begin{gather*}
 \bar\nabla_X U = \# (\nabla_{\#X} Y + \{ X,Y\}),
 \end{gather*}
 where $\{X,Y\}:= \nabla_{\#Y}X- \nabla_{\#X}Y+[X,Y]$. Evidently this representation pulls back, under $f$, to the representation of $A(f)$ on $T_\Sigma M \cong E/{\mathfrak h} $ induced by the fundamental representation~\eqref{hsq}.
\end{proof}

Globalisations of $A(f)^j$ for $j \ge 3$ are not described here.

\subsection*{Intrinsic geometry}
Suppose $j\ge 1$ and that the Lie algebroid $A^j$ in an invariant filtration is transitive. Then, according to Proposition~\ref{cP}, $A^{j+1}$ acts canonically on $T \Sigma $. Such a representation equips $\Sigma $ with an `intrinsic infinitesimal geometry' in the following way. View the representation as a Lie algebroid morphism $A^{j+1} \rightarrow \mathfrak{gl}[T \Sigma] $, where $\mathfrak{gl}[T \Sigma]$ is the Lie algebroid of derivations on $T \Sigma $ (the Lie algebroid of the frame groupoid of $T \Sigma $). The adjoint representation $\adjoint \colon J^1 (T \Sigma ) \rightarrow \mathfrak{gl}[T \Sigma ] $, defined by $\adjoint_{J^1 X}Y=[X,Y]$, is an isomorphism, allowing us to regard the image of $A^{j+1} \rightarrow \mathfrak{gl}[T \Sigma] $ as a subalgebroid of $J^1 (T \Sigma)$, i.e., as in infinitesimal geometric structure, in the sense of~\cite{Blaom_12}. For example, if $A^{j+1} \cong T \Sigma $ then its representation on $T \Sigma $ amounts to an infinitesimal parallelism on $\Sigma $. In the case of a~hypersurface $f \colon \Sigma \hookrightarrow M$ in a Riemannian space form, the representation of $A^2(f)$ on $T \Sigma $ determines a~representation on $S^2\,T^*\Sigma$ whose sole invariant sections are constant multiples of the inherited metric (see Remark~\ref{ifR}).

\subsection*{Planar curves}
For a simple illustration of the theory now developed, let us recover the following well-known fact: A unit-speed planar curve $f \colon I \rightarrow {\mathbb R}^2$, defined on some interval $I$, is completely characterized, up to orientation-preserving rigid motions, by its curvature $\kappa $, an invariant under such motions.

Let $\xi$, $\eta $ denote the standard coordinate functions on ${\mathbb R}^2$ and $x,y \colon I \rightarrow {\mathbb R}^2$ their pullbacks under a regular curve $f \colon I \rightarrow {\mathbb R}^2 $. The plane is a Klein geometry with transitively acting group~$G$ the group of orientation-preserving rigid motions, motions that preserve the speed of curves. The corresponding Lie algebra ${\mathfrak g} $ of planar vector fields admits the basis $\{e_1,e_2,e_3\}$, where
\begin{gather*}
 e_1 =\frac{\partial }{\partial \xi },\qquad e_2= \frac{\partial }{\partial \eta},\qquad e_3= -\eta \frac{\partial }{\partial \xi }+\xi \frac{\partial }{\partial \eta}.
\end{gather*}
We have $[e_1,e_2]=0$, $[e_3,e_1]=-e_2$, and $[e_3,e_2]=e_1$.

Anticipating that polynomial invariants on ${\mathfrak g} $ are going to play a role, we let $Q$ denote the positive definite quadratic form on the commutator $[{\mathfrak g}, {\mathfrak g}]=\spann\{e_1,e_2\}$, with respect to which~$e_1$ and~$e_2$ form an orthonormal basis, and let $\mu \colon {\mathfrak g} \rightarrow {\mathbb R} $ be the linear form vanishing on $[{\mathfrak g},{\mathfrak g}]$ with $\mu(e_3)=1$. Both $Q$ and $\mu $ are $\Adjoint$-invariant, and these invariants encode the structure of ${\mathfrak g} $: Agreeing to give $[{\mathfrak g}, {\mathfrak g}]$ the orientation defined by the ordered basis $(e_1,e_2)$, then for any other positively-oriented, $Q$-orthogonal basis $({\mathbf x},{\mathbf y})$, where ${\mathbf x}$ and ${\mathbf y}$ have the same $Q$-length, we have $[{\mathbf x},{\mathbf y}]=0$ and, for all ${\mathbf z} \in {\mathfrak g} $,
\begin{gather}
 [{\mathbf z}, {\mathbf x}]= - \mu({\mathbf z}){\mathbf y}, \qquad [{\mathbf z}, {\mathbf y}] = \mu({\mathbf z}) {\mathbf x}.\label{sense}
\end{gather}
Viewing some ${\mathbf z} \in {\mathfrak g} $ concretely as a Killing field, $\mu({\mathbf z})$ is the curl of ${\mathbf z}$. Of course $[{\mathfrak g}, {\mathfrak g} ] $ is just the space of constant vector fields on ${\mathbb R}^2 $. Bearing Remark~\ref{polynomialR} in mind, we record for later use:
\begin{Lemma}\label{lemmaA}
 Given $\lambda >0$ there exists an automorphism of ${\mathfrak g} $ pulling $Q$ back to $\lambda Q$.
\end{Lemma}
\begin{proof} Consider the automorphism $e_1 \mapsto \lambda e_1$, $e_2 \mapsto \lambda e_2$, $e_3 \mapsto e_3$.
\end{proof}

Since $f $ is regular (an immersion) we may identify $A(f)$ with the subbundle of $E:={\mathfrak g} \times I$ consisting of those $(\zeta,t)$ for which the Killing field $\zeta(t)$ is tangent to $f$ at $f(t)$. If
$\{E,E\}:=[{\mathfrak g},{\mathfrak g}]\times I \subset E$, then the invariants $Q$ and $\mu $ define ${\mathbb R} $-valued maps on $\{E,E\}$ and $E$, which we denote by the same symbols. Evidently, $A(f)$ has rank two and there exists a unique section $X$ of $A(f)\cap \{E,E\}$ such that $X$ has constant $Q$-length one, and such that $\# X = \frac{1}{s} \frac{\partial}{\partial t}$ for some positive function $s \colon I \rightarrow (0,\infty)$, namely the speed of $f$:
\begin{gather*}
 X = \frac{\dot x}{s} \frac{\partial }{\partial \xi } + \frac{\dot y}{s} \frac{\partial }{\partial \eta}.
\end{gather*}

Next, choose $Y \in \Gamma(\{E,E\})$ such that $(X,Y)$ forms a positively oriented $Q$-orthonormal basis of each fibre of $\{E,E\}$:
\begin{gather*}
 Y = -\frac{\dot y}{s} \frac{\partial }{\partial \xi } + \frac{\dot x}{s} \frac{\partial }{\partial \eta}.
\end{gather*}
The map $Q \colon \{E,E\}\rightarrow {\mathbb R} $ is tautologically invariant with respect to the canonical flat connec\-tion~$\nabla $. We now differentiate $Q(X,X)=1$ to obtain $Q(X, \nabla_{\partial/\partial t}X)=0$, implying that
\begin{gather}
 \nabla_{\partial /\partial t}X= s\kappa Y,\label{jert1}
\end{gather}
for some function $\kappa \colon I \rightarrow {\mathbb R} $, namely the curvature. Differentiating the identities $Q(X,Y)=0 $ and $Q(Y,Y)=1$, we deduce that
\begin{gather}
 \nabla_{\partial /\partial t}Y= -s\kappa X.\label{jert2}
\end{gather}
Finally, a section $Z$ of the rank-one kernel of the anchor $\# \colon A(f) \rightarrow TI$ is fixed by requiring $\mu(Z)=1$. Since the Killing field $Z(t)$ must vanish at $f(t)=(x(t),y(t))$, we obtain
\begin{gather*}
 Z=(y- \eta)\frac{\partial }{\partial \xi }-(x-\xi) \frac{\partial
 }{\partial \eta}. 
\end{gather*}
Since $\partial/\partial \xi $ and $\partial/\partial \eta$ are constant sections of ${\mathfrak g} \times I$, we obtain, on differentiating,
\begin{gather}
 \nabla_{\partial/\partial t} Z = \dot y \frac{\partial }{\partial \eta }- \dot x \frac{\partial }{\partial \eta }=-sY. \label{jert3}
\end{gather}

The sections $X$, $Y$, $Z$ furnish a basis for each fibre of $E={\mathfrak g} \times I$ and $A(f)$ is spanned by $X$ and $Z$. The algebraic bracket on $E$ is given by
\begin{gather}
 \{X,Y\}=0, \qquad \{Z,X\}=-Y, \qquad \{Z,Y\}=X,\label{cew}
\end{gather}
where the last two relations follow from \eqref{sense} or direct computation. Since $\#X=\frac{1}{s}\frac{\partial }{\partial t}$ and \smash{$\#Z=0$}, the Maurer--Cartan equations \eqref{mce} show that $[X,Z]=0$, which completes the description of the Lie algebroid structure of $A(f)$.

\begin{figure}[t]\centering
 \includegraphics[scale=0.3]{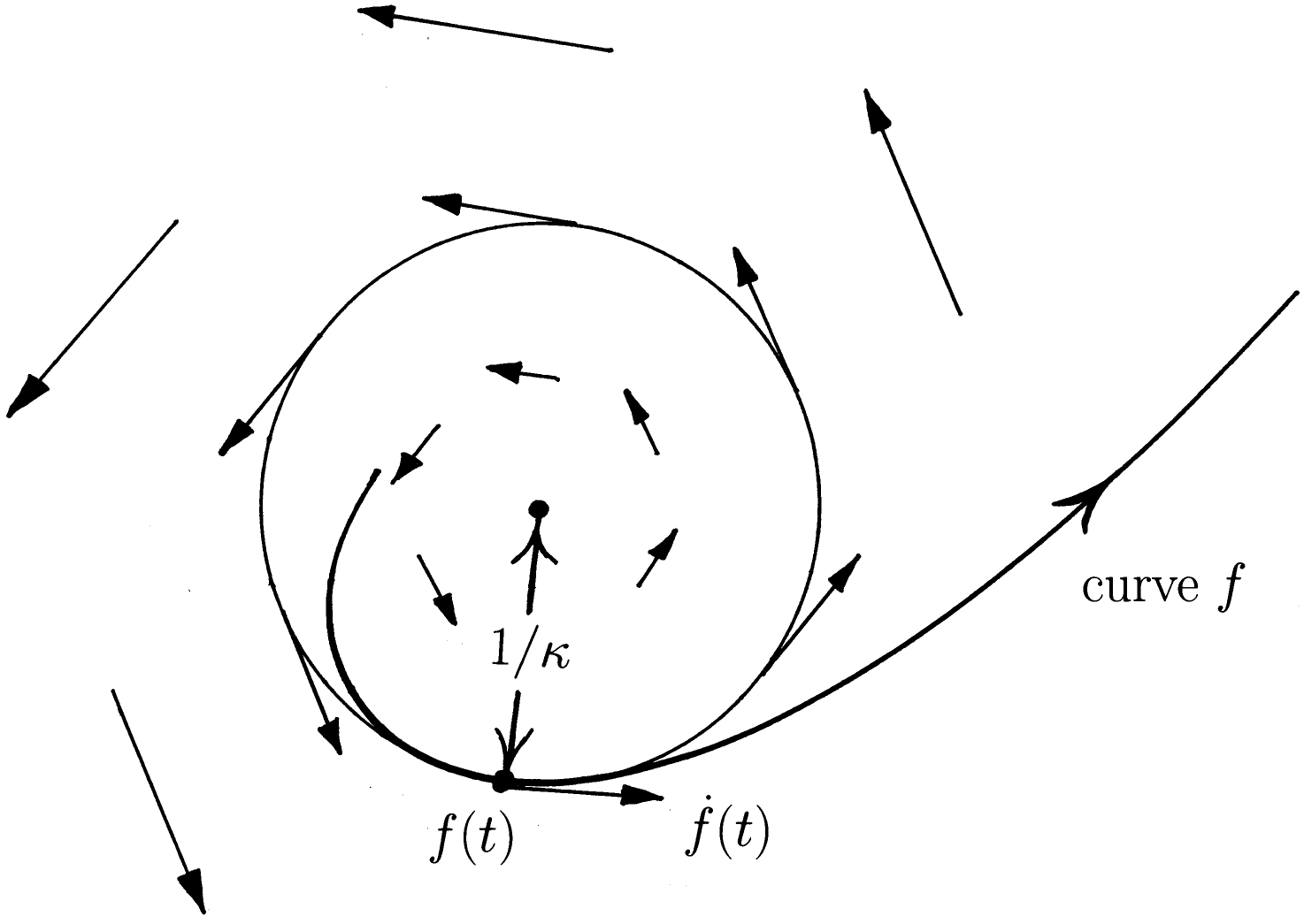}
\caption{The vector field $X^2(t)\in {\mathfrak g}$, where $X^2=s(X + \kappa Z$). }\label{Fig1}
\end{figure}

From \eqref{jert1}--\eqref{jert3}, and our definition of the invariant filtration, we see that $A(f)^2$ is generated by the single section $X^2:=s(X + \kappa Z)$, which the anchor maps to $\partial /\partial t$, giving us an identification $A(f)^2 \cong TI$.
\begin{Remark}\label{planarR} Notice that $X^2$ has the following geometric interpretation: Evaluated at $t \in I$, it is the planar Killing field whose integral curve through $f(t)$ is the circle best approximating the
 curve $f$ at $f(t)$ (Fig.~\ref{Fig1}).
\end{Remark}

One computes $\nabla_{\partial /\partial t}X^2= \dot s X + (s \kappa)\dot{~} Z$, from which it follows that $A(f)^3=A(f)^2$ when the speed and curvature are constant, and $A(f)^3=0$ otherwise. So, according to Theorem~\ref{sssT}, we recover unit speed circular arcs and straight line segments as the curves with non-trivial symmetry.

With the help of \eqref{jert1}--\eqref{cew} we can compute the fundamental representation; it is given by
\begin{alignat*}{4}
& \bar\nabla_XX =\kappa Y,\qquad&& \bar\nabla_XY =-\kappa X,\qquad&& \bar\nabla_XZ= 0,&\\
& \bar\nabla_ZX =-Y, \qquad & & \bar\nabla_ZY =X, \qquad & & \bar\nabla_Z Z=0.&
\end{alignat*}
We leave the reader to show that the canonical representation of $A(f)^2 \cong {T}{I}$ on $TI$ described in Proposition \ref{cP} is given by
\[
\bar\nabla_{\partial /\partial t}\,\frac{\partial }{ \partial t} =\frac{\dot s}{s}\frac{\partial }{\partial t}.
\]
In particular, $\partial/ \partial t$ is parallel in the intrinsic geometry inherited by $I$ precisely when the speed is constant. (Intrinsic geometry is discussed in Section~\ref{invint}.)

Our next task is to construct invariants of arbitrary infinitesimal immersions $\omega \colon A \rightarrow {\mathfrak g} $. To this end, we state:
\begin{Lemma}\label{lemmaB} Let $\omega \colon A \rightarrow {\mathfrak g} $ be any infinitesimal immersion of $I $ into ${\mathbb R}^2$ and use $\omega $ to identify $A$ with a subbundle of $E={\mathfrak g} \times I$. Then $A \cap \{E,E\}$ has rank one and is mapped by the anchor $\# \colon A \rightarrow TI$ onto $TI$.
\end{Lemma}
\begin{proof} Since $I$ is simply-connected, $\omega $ admits a primitive $f \colon I \rightarrow {\mathbb R}^2$, by Theorem~\ref{infirmT}. Since $\omega $ and $\delta f$ are isomorphic, the lemma follows because it is already true for the logarithmic derivatives of curves.
\end{proof}

For an infinitesimal immersion $\omega \colon A \rightarrow {\mathfrak g} $ we can accordingly find a~unique section $X$ of $A \cap \{E,E\}$ with constant $Q$-length one such that $\# X= \frac{1}{s}\frac{\partial }{\partial t}$ for some positive function $s \colon I \rightarrow (0,\infty)$ we call the {\df speed} of $\omega $. We choose a second section $Y$ of $\{E,E\}$ as we did for curves and, arguing as before, obtain a function $\kappa \colon I \rightarrow {\mathbb R} $ defined by~\eqref{jert1} that we call the {\df curvature} of $\omega $. By construction, the speed and curvature of the logarithmic derivative $\delta f$ of a~curve $f$ coincide with the usual speed and curvature of $f$. It is clear that speed and curvature are invariants of an infinitesimal immersion.

We remark that the curvature of $\omega $ can be equivalently defined by $\kappa =\mu\big(X^2\big)/s$ where $X^2$ is the unique section of $A^2$ such that $\#X^2=\partial /\partial t$.

We now apply the Construction Principle \ref{constructionPrinciple} to obtain:
\begin{Theorem}\label{planarT} For every smooth function $\kappa \colon I \rightarrow {\mathbb R}$ there exists a unit-speed curve $f \colon I \rightarrow {\mathbb R}^2 $ with curvature $\kappa $, unique up to orientation-preserving isometries.
\end{Theorem}
\begin{proof} Let $X$, $Y$, $Z$ denote the constant sections of the trivial bundle $E := {\mathbb R}^3 \times I $ and define an algebraic bracket $\{\,\cdot\,,\,\cdot\,\}$ on $E$ by the relations \eqref{cew}. This makes $E $ into a ${\mathfrak g} $-bundle, with~${\mathfrak g} $ defined as above. Define a connection $\nabla $ on $E$ by \eqref{jert1}--\eqref{jert3}, with $s=1$, and verify that the algebraic bracket is $\nabla $-parallel. Make the subbundle $A \subset E$ spanned by $X$ and $Z$ into a~Lie algebroid by declaring $[X,Z]=0$, $\#X= \partial /\partial t$ and $\#Z=0$. Then $A$ is integrable and the Maurer--Cartan equations~\eqref{mce} hold. Applying Proposition~\ref{croonP}, we obtain an isomorphism $E \cong {\mathfrak g} \times I$ such that the composite of $A \hookrightarrow E $ with the projection $E \rightarrow {\mathfrak g} $ is Lie algebroid morphism $\omega \colon A \rightarrow {\mathfrak g} $. It is easy to see that $\omega \colon A \rightarrow {\mathfrak g} $ must be an infinitesimal immersion. Since $\{E,E\} \subset E$ must be $\nabla $-invariant, we have $\nabla_{\partial / \partial t}X \in \Gamma(\{E,E\})$. Since $Q$, viewed as an inner product on $E \cong {\mathfrak g} \times I$, is $\nabla $-parallel, we compute
 \begin{gather*}
 \frac{\partial }{\partial t}Q(X,X)=2Q(X,\nabla_{\partial/\partial t} X)=2 \kappa\, Q(X,Y).
 \end{gather*}
 On the other hand, since $Q$ is symmetric and $\adjoint$-invariant, we have
 \begin{gather*}
 0=Q(X,\{Z,X\})= Q(X,Y).
 \end{gather*}
Combining these two equations, we see that $Q(X,X)$ is constant. Applying Lemma \ref{lemmaA}, we may arrange, by composing $\omega \colon A \rightarrow {\mathfrak g} $ with an outer automorphism of ${\mathfrak g} $ if necessary, that $Q(X,X)=1$.

By construction, $\omega $ has speed one and curvature $\kappa $. By Theorem \ref{infirmT}, $\omega $ has a primitive $f \colon I \rightarrow {\mathbb R}^2$, unique up to orientation-preserving isometries of ${\mathbb R}^2$. Since $\omega $ and $\delta f$ are isomorphic, they must have the same invariants, so that $f$ is a unit-speed curve with curvature $\kappa $.
\end{proof}

\subsection*{Curves in the elliptic plane}
We now sketch a similar analysis of curves on $M=S^2$ with $G=\operatorname{SO}(3)$. In this case all polynomial invariants are generated by the Killing form of ${\mathfrak g} $, which we scale by a factor of $-1/2$ to obtain an inner product $Q$ on ${\mathfrak g} $, and on the trivial bundle $E={\mathfrak g} \times I$.

Given an infinitesimal immersion $\omega \colon A \rightarrow {\mathfrak g}$ we define a section $X$ of the rank-two Lie algebroid $A$ by requiring it to be $Q$-orthogonal to the kernel ${\mathfrak h} $ of the anchor, to have constant $Q$-length one, and satisfy $\# X = \frac{1}{s}\frac{\partial }{ \partial t}$, for some positive function $s \colon I \rightarrow {\mathbb R} $ called the {\df speed} of~$\omega $. This speed coincides with the
usual speed of a curve $f \colon I \rightarrow S^2$ when we take~$\omega $ to be the logarithmic derivative of a curve $f$ (with our choice of scaling factor, $-1/2$). Let $Y$ be a~section of~$E$ that is $Q$-orthogonal to $A$ and has constant $Q$-length one. Such a $Y$ is determined uniquely up to sign. Differentiating $Q(Y,Y)=1$, one shows that $\nabla_{\partial/\partial t}Y$ is a section of $A$ and can therefore define $\kappa \colon I \rightarrow {\mathbb R} $ by
\begin{gather*}
 \#\nabla_{\partial/\partial t}Y=-\kappa\frac{\partial }{\partial t}.
\end{gather*}
This, in turn, defines the {\df curvature} $|\kappa|$ of $\omega $, which is independent of the choice of $Y$ above, and hence an invariant, and which coincides with the absolute value of the curvature of a curve $f \colon I \rightarrow S^2$ when we take $\omega $ to be the logarithmic derivative of $f$. Note that the usual (signed) curvature of a curve is not an invariant in our theory because symmetries of $S^2$ include the orientation-reversing isometries.

For a curve $f$ one goes on to define $Z=\{X,Y\}$, a section of ${\mathfrak h} $, and then obtains a complete description of the logarithmic derivative of $f$ readily enough. One then proves an analogue of Theorem~\ref{planarT} by applying Construction Principle~\ref{constructionPrinciple}. Details are left to the reader.

For an analysis of curves in the equi-affine plane, see Section~\ref{appendixA}.

\section{The invariant filtration for Riemannian subgeometry}\label{bonnetsec2}
In this section we study oriented codimension-one immersed submanifolds $\Sigma $ of the Riemannian manifold $M = S^n$, ${\mathbb R}^n$, or ${\mathbb H}^n $.

\subsection*{The Lie algebroid associated with a Riemannian manifold}
The basic groupoid associated with an oriented Riemannian manifold $M$ is the Lie groupoid $\mathrm{SO}[TM]$ of orientation-preserving orthonormal relative frames. Being a frame groupoid, the Lie algebroid of~$\mathrm{SO}[TM]$ can be viewed as a~subalgebroid $A $ of the Lie algebroid $\mathfrak{gl}[TM]$ of derivations of $TM$ (defined in, e.g., \cite[Section~3.3]{Mackenzie_05}). However, the Levi-Cevita connection~$\nabla^M$ furnishes a splitting $u \mapsto \nabla_u^M$ of the canonical exact sequence
\begin{gather*}
 0 \rightarrow \mathfrak{so}(TM) \rightarrow A \xrightarrow{\#} TM \rightarrow 0,
\end{gather*}
allowing us to identify $A$ with
\begin{gather*}
 \mathfrak{so}[TM]:= {TM} \oplus \mathfrak{so}({TM}).
\end{gather*}
Recall here that $\mathfrak{so}({TM}) \subset T^*\!M \otimes TM$ denotes the subbundle of tangent space endomorphisms that are skew-symmetric with respect to the metric -- a Lie algebra bundle of type $\mathfrak{so}(n)$, $n=\dim M$.

The anchor of $\mathfrak{so}[TM]$ is projection onto the first summand. Writing elements and sections of $\mathfrak{so}[TM]={TM} \oplus \mathfrak{so}({TM})$ vertically, the Lie bracket is given by
\begin{gather*}
 \left[
 \begin{pmatrix}
 V_1\\ \Phi_1
 \end{pmatrix},
 \begin{pmatrix}
 V_2\\ \Phi_2
 \end{pmatrix}
\right] =
\begin{pmatrix}
 [V_1,V_2]\\
 [\Phi_1, \Phi_2] + \nabla^M_{V_1} \Phi_2 - \nabla^M_{V_2}
 \Phi_1 + \curvature \nabla^M(V_1,V_2)
\end{pmatrix},
\end{gather*}
as the reader is invited to check. The canonical representation of $\mathfrak{so}[TM]$ on $TM$ is given in the present model by
\begin{gather}
 \bar\nabla_{
 \left(\begin{smallmatrix}
 U \\ \Phi
 \end{smallmatrix}\right)}V=\nabla^M_U V + \Phi V.\label{brook}
\end{gather}
This representation induces a representation of $\mathfrak{so}[TM]$ on the bundle $S^2\, T^*\!M$ of which the metric on $M$ is a section. By construction, the metric is invariant under this representation and is, up to scale, the unique such section.

As explained in Section \ref{introduction}, we are also interested in the Lie groupoid $\mathrm{SO}[T^+M]$ of orientation-preserving orthonormal relative frames of $T^+M:=TM \oplus ({\mathbb R} \times M)$, equipped with the inner product
\begin{gather}
 \langle\!\langle V_1 \oplus a_1,V_2 \oplus a_2\rangle\!\rangle := \langle\!\langle V_1,V_2\rangle\!\rangle + a_1 a_2.\label{cds}
\end{gather}
This product is $\nabla^{M+}$-parallel if we define $\nabla^{M+}_u(V \oplus a)=\nabla^M_u V \oplus da(u)$ and, by a~similar argument, we obtain a split exact sequence for the Lie algebroid of $\mathrm{SO}[T^+M]$, identifying it with $TM \oplus \mathfrak{so}(T^+ M)$, which in turn has a natural identification with
\[
\mathfrak{so}[T^+M]:=TM \oplus \mathfrak{so}(TM) \oplus TM =\mathfrak{so}[TM]\oplus TM.
\]
The anchor of $\mathfrak{so}[T^+M]$ is given by $\#(X, w)=\#X$ and the bracket by
\begin{gather*}
 [X_1 \oplus W_1, X_2 \oplus W_2] = \left([X_1,X_2] +
 \begin{pmatrix}
 0 \\ W_1^\flat \otimes W_2 - W_2^\flat \otimes W_1
 \end{pmatrix}\right)\oplus \big(\bar\nabla_{X_1} W_2 -
 \bar\nabla_{X_2} W_1\big),
\end{gather*}
where $\bar\nabla $ denotes the canonical representation of $\mathfrak{so}[TM]$ on $TM$ defined in~\eqref{brook}. The Lie algebroid $\mathfrak{so}[T^+M]$ itself acts on $T^+M $ according to
\begin{gather*}
 \bar\nabla_{X \oplus W}(V \oplus a)= (\bar\nabla_X V + aW)\oplus ({\rm d}a(\#X)- \langle\!\langle W,V\rangle\!\rangle).
\end{gather*}
By construction, this action preserves the inner product \eqref{cds} which, up to scale, is the unique such product.

\subsection*{The Cartan connection for Riemannian geometry}
There is a canonical Cartan connection $\nabla $ on $\mathfrak{so}[TM]$ \cite{Blaom_06,Blaom_12}, given by
\begin{gather}
 \nabla_U \begin{pmatrix}
 V \\ \Phi
 \end{pmatrix} =
 \begin{pmatrix}
 \nabla^M_UV + \Phi U \\
 \nabla^M_U\Phi + \curvature \nabla^M (U,V)
 \end{pmatrix}.\label{form}
\end{gather}
Since $\nabla $ is a Cartan connection, the space ${\mathfrak g}_\nabla$, consisting of all $\nabla $-parallel sections of $\mathfrak{so}[TM]$, forms a Lie subalgebra of all sections.

To recall the significance of $\nabla $, remember that a vector field $\xi $ on $M$ is a Killing field (infinitesimal isometry) precisely when the corresponding section $\nabla^M \xi $ of $T^*\!M \otimes {TM} $ is in fact a~section of $\mathfrak{so}({TM})$. The following is proven in~\cite{Blaom_12}:
\begin{Theorem}\label{rias2T}
 For every Killing field $\xi $ the section
 \[
 X_\xi = \begin{pmatrix}
 \xi \\ -\nabla^M \xi
 \end{pmatrix}
 \]
 of $\mathfrak{so}[TM]$ is $\nabla $-parallel and all parallel sections are of this form. In particular, the anchor $\# \colon \mathfrak{so}[TM] \rightarrow TM$, when viewed as a map of corresponding section spaces $\Gamma(\mathfrak{so}[TM]) \rightarrow \Gamma({TM})$, maps the Lie algebra ${\mathfrak g}_\nabla$ of $\nabla $-parallel sections of $\mathfrak{so}[TM]$ isomorphically onto the Lie algebra of Killing fields of~$M$.
\end{Theorem}

By this result the curvature of $\nabla $ (see below) is the local obstruction to the existence of Killing fields.

The bracket on ${\mathfrak g}_\nabla$ can be expressed in terms of the {\df torsion} of $\nabla $ -- an algebraic bracket defined on any Lie algebroid equipped with a connection by
\[\{X,Y\}:=\nabla_{\#Y}X - \nabla_{\#X}Y + [X,Y].\]
For if $\xi,\eta \in {\mathfrak g}_\nabla $ then $[\xi,\eta]=\{\xi,\eta\}$. In the present case one computes
\begin{gather}
 \left\{
 \begin{pmatrix}
 V_1 \\ \Phi_1
 \end{pmatrix},
 \begin{pmatrix}
 V_2 \\ \Phi_2
 \end{pmatrix}
 \right\}=
 \begin{pmatrix}
 \Phi_1 V_2 -\Phi_2 V_1 \\
 [\Phi_1,\Phi_2] - \curvature \nabla^M (V_1,V_2)
 \end{pmatrix}. \label{torsion}
\end{gather}

The curvature of the connection $\nabla $ on $\mathfrak{so}[TM]$ is readily computed and seen to vanish precisely when $M$ has constant sectional curvature~$s$, i.e., when the curvature of the Levi-Cevita connection is given by
\begin{gather}
 \curvature \nabla^M (V_1,V_2)=s\big(V_2^\flat \otimes V_1 - V_1^\flat \otimes V_2\big).\label{torsion2}
\end{gather}
Here $V \mapsto V^\flat \colon TM \rightarrow T^*\!M$ is the canonical isomorphism determined by the metric. Assu\-ming~$M$ is simply-connected, the Lie algebra ${\mathfrak g}_\nabla $ assumes its maximal possible dimension $k$ in this case; $k=\rank \mathfrak{so}[TM]=n(n+1)/2 $, where $n=\dim M$.

\subsection*{The ambient space of embeddings}
In the remainder of this section $M$ denotes either the sphere $S^n$, Euclidean space ${\mathbb R}^n $, or hyperbolic space ${\mathbb H}^n$. Then $M$ is a Klein geometry with transitively acting group $G$, where $G$ is $\SO(n+1)$, the semidirect product $\SO(n) \ltimes {\mathbb R}^n $, or $\SO^+(n,1)$, respectively. Note that $M$ has constant sectional curvature $s=1,0,-1$, respectively. According to \cite[Examples~2.3]{Blaom_F}, $\automorphism(M)=G$ in all cases except the even-dimensional spheres, in which case $\automorphism(M)$ is the group of all isometries, both orientation-preserving and orientation-reversing. The Lie algebra of $G$ will be identified with the Lie algebra ${\mathfrak g} \subset \Gamma(TM)$ of all Killing fields (infinitesimal isometries), which is possible because $G$ acts faithfully and $\dim(G)=n(n+1)/2$.
\begin{Proposition}\label{ambientP} The action algebroid ${\mathfrak g} \times M$ is isomorphic to the Lie algebroid $\mathfrak{so}[TM]$ defined above. An explicit isomorphism is given by
 \[ (\xi, m) \mapsto X_\xi (m):= \begin{pmatrix}
 \xi(m) \\ -(\nabla^M \xi)(m)
 \end{pmatrix}\colon \ {\mathfrak g} \times M \rightarrow \mathfrak{so}[TM].
 \]%
 Under this isomorphism the canonical flat connection on ${\mathfrak g} \times M$ coincides with the Cartan connection $\nabla $.
\end{Proposition}

\begin{proof} It is readily checked that the map is a Lie algebroid morphism. To show the map is an isomorphism it suffices, by a dimension count, to establish surjectivity. But every element of $\mathfrak{so}[TM]$ can be extended to a $\nabla $-parallel section, because $\nabla $ is flat ($M$ has constant scalar curvature) and $M$ is simply-connected. By the preceding theorem this section is of the form $X_\xi$ for some $\xi \in {\mathfrak g}$. The last statement in the theorem is now obvious.
\end{proof}

\subsection*{The logarithmic derivative of an immersion}
Now let $\Sigma$ denote an immersed, oriented, codimension-one submanifold, and $f \colon \Sigma \rightarrow M$ the immersion. The logarithmic derivative of $f$ is a map $\delta f \colon A(f) \rightarrow {\mathfrak g} $, where $A(f)$ is the pullback of the action algebroid ${\mathfrak g} \times M$ by $f$. In view of the preceding proposition, we may identity ${\mathfrak g} \times M $ with $\mathfrak{so}[TM]={TM} \oplus \mathfrak{so}({TM})$. Then, as the anchor of $\mathfrak{so}[TM]$ is just projection onto the first summand, we have
\begin{gather*}
 A(f)=T \Sigma \oplus \mathfrak{so}({T_\Sigma M}),
\end{gather*}
where $T_\Sigma M $ denotes the vector bundle pullback of ${TM} $, and $\mathfrak{so}({T_\Sigma M})$ is the bundle of skew-symmetric endomorphisms thereof. In the next subsection we shall show that the Lie bracket on $A(f)$ is given by
\begin{gather}
 \left[\begin{pmatrix} v_1 \\ \Phi_1\end{pmatrix}, \begin{pmatrix}
 v_2 \\ \Phi_2\end{pmatrix}\right]=
 \begin{pmatrix}
 [v_1, v_2] \\
 [\Phi_1,\Phi_2] + \nabla^M_{v_1} \Phi_2 - \nabla^M_{v_2} \Phi_1 + s\big(v_2^\flat \otimes v_1 - v_1^\flat \otimes v_2\big)
 \end{pmatrix}.\label{logs:1}
\end{gather}
In this case we should interpret $\nabla ^M$ as the pullback of the Levi-Cevita connection on $TM$ to a connection on $T_\Sigma M$ (and a~corresponding connection on $\mathfrak{so}(T_\Sigma M)$). We use lowercase letters to distinguish elements of $T \Sigma $ from elements of $TM$. Under our identifications, the logarithmic derivative $\delta f$ of $f$ is just the composite of the inclusion
\[T \Sigma \oplus \mathfrak{so}(T_\Sigma M) \monomorphism {TM} \oplus \mathfrak{so}({TM})=\mathfrak{so}[TM]\] with the projection $\mathfrak{so}[TM]\cong {\mathfrak g} \times M \rightarrow {\mathfrak g} $.

\subsection*{The fundamental representation}
As in earlier sections, we now regard $A(f)$ as a subbundle of $E={\mathfrak g} \times \Sigma $. As ${\mathfrak g} \times \Sigma $ is just the vector bundle pullback of ${\mathfrak g} \times M \cong \mathfrak{so}[TM]=TM \oplus \mathfrak{so}(TM) $, we may make the identification
\[E =T_\Sigma M \oplus \mathfrak{so}(T_\Sigma M).\] Then the inclusion of $A(f)=T \Sigma \oplus \mathfrak{so}(T_\Sigma M)$ into $E$ is the obvious one.

According to Proposition \ref{ambientP}, the canonical flat connection on ${\mathfrak g} \times M \cong \mathfrak{so}[TM]$ is represented by the Cartan connection $\nabla $ on $\mathfrak{so}[TM]$ defined by~\eqref{form}, and this connection pulls back to a flat connection on $E=T_\Sigma M \oplus \mathfrak{so}(T_\Sigma M)$ given by the same formula, provided we again interpret~$\nabla^M$ as the pullback of the Levi-Civita connection to a connection on $T_\Sigma M$.

Under the identification ${\mathfrak g} \times M \cong \mathfrak{so}[TM]$ the canonical algebraic bracket on sections of ${\mathfrak g} \times M$ coincides with the torsion of $\nabla$, which is given explicitly in~\eqref{torsion}, \eqref{torsion2} above. From this we may determine the corresponding algebraic bracket on sections of~$E$:
\begin{gather}
 \left\{
 \begin{pmatrix}
 V_1 \\ \Phi_1
 \end{pmatrix},
 \begin{pmatrix}
 V_2 \\ \Phi_2
 \end{pmatrix}
 \right\}=
 \begin{pmatrix}
 \Phi_1 V_2 -\Phi_2 V_1 \\
 [\Phi_1,\Phi_2] + s\big(V_1^\flat \otimes V_2 - V_2^\flat \otimes V_2\big)
 \end{pmatrix}. \label{torsion3}
\end{gather}
This bracket is necessarily $\nabla $-invariant. Since, furthermore, the Maurer--Cartan equations \eqref{mce} must hold, one now readily establishes the earlier claim \eqref{logs:1} regarding the Lie algebroid bracket on $A(f)$.

With the help of formula \eqref{form} for the connection $\nabla $ and formula \eqref{torsion3} for the algebraic bracket on~$E$, the fundamental representation of $A(f)=T \Sigma \oplus \mathfrak{so}(T_\Sigma M)$ on $E=T_\Sigma M \oplus \mathfrak{so}(T_\Sigma M)$ is readily computed and given by
\begin{gather}
 \bar\nabla_{
\left( \begin{smallmatrix}
 u \\ \Phi
 \end{smallmatrix}\right)}
 \begin{pmatrix}
 V \\ \Psi
 \end{pmatrix}=
 \begin{pmatrix}
 \nabla^M_u V + \Phi V \\ \nabla^M_u \Psi + [\Phi,\Psi]
 \end{pmatrix}.\label{seven}
\end{gather}

\subsection*{Decomposing $\boldsymbol{T_\Sigma M}$}
Since $M$ and $\Sigma $ are oriented, $\Sigma$ has a well-defined unit normal vector field, leading to identifications
\begin{gather}
 T_\Sigma M \cong T \Sigma \oplus ({\mathbb R} \times \Sigma), \qquad \mathfrak{so}(T_\Sigma M) \cong \mathfrak{so} (T \Sigma) \oplus T \Sigma.\label{wpo}
\end{gather}
Elements of the bundles above will be written horizontally. Under the above identifications the tautological representation of $\mathfrak{so}(T_\Sigma M)$ on $T_\Sigma M$ becomes a representation of~$\mathfrak{so}(T\Sigma) \oplus T\Sigma $ on $T\Sigma \oplus ({\mathbb R} \times \Sigma)$, here denoted by~$\adjoint$, and given by
\begin{gather*}
 \adjoint_{\phi \oplus w}(v \oplus a)=(\phi v + a w)\oplus - \langle\!\langle w,v\rangle\!\rangle,
\end{gather*}
where $\langle\!\langle \, \cdot \,,\,\cdot \,\rangle\!\rangle$ denotes the metric. The bracket on $\mathfrak{so}(T_\Sigma M)$ becomes a bracket on $\mathfrak{so}(T \Sigma) \oplus T \Sigma$ given by
\begin{gather*}
 [\phi_1 \oplus w_1, \phi_2 \oplus w_2]=\big([\phi_1,\phi_2] + w_1^\flat \otimes w_2 - w_2^\flat \otimes w_1\big)\oplus (\phi_1 w_2 - \phi_2 w_1).
\end{gather*}
For the reader's convenience we also record well-known expressions for the pullback of the Levi-Cevita connection $\nabla^M$ to connections on $T_\Sigma M$ and $\mathfrak{so} (T_\Sigma M)$, after making the identifications in~\eqref{wpo}:
\begin{gather*}
 \nabla^M_u (v \oplus a) =\big(\nabla_u^\Sigma v - a \ii(u)^\sharp\big)\oplus ({\rm d}a(u) + \ii(u,v)),\\
 \nabla^M_u (\phi \oplus w) = \big(\nabla^\Sigma_u \phi + w^\flat \otimes \ii(u)^\sharp - \ii(u)\otimes w\big) \oplus \big( \nabla^\Sigma_u w + \phi\big(\ii(u)^\sharp\big)\big).
\end{gather*}
In these formulas $\nabla^\Sigma $ denotes the Levi-Cevita connection on $\Sigma $ and $\ii $ the second fundamental form. The map $V \mapsto V^\flat \colon TM \rightarrow T^*\!M$ is the canonical isomorphism determined by the metric, and $\alpha \mapsto \alpha ^\sharp $ is its inverse.

Having made the identifications in~\eqref{wpo} we have the following updated model of~$E$:
\begin{gather} E=
\begin{matrix}
 T \Sigma \oplus ({\mathbb R} \times M) \\ \oplus \\ \mathfrak{so}(T \Sigma) \oplus T \Sigma. \end{matrix} \label{mmmm}
\end{gather}
We are now writing bundle direct sums as we will be laying out corresponding elements. The algebraic Lie bracket on $E$ is now given by the formula
\begin{gather}
 \left\{
 \begin{pmatrix}
 v_1 \oplus a_1 \\ \phi_1 \oplus w_1
 \end{pmatrix},
 \begin{pmatrix}
 v_2 \oplus a_2 \\ \phi_2 \oplus w_2
 \end{pmatrix}
 \right\}=
 \begin{pmatrix}
 \text{\circled{1}}\oplus\text{\circled{2}} \\ \text{\circled{3}}
 \oplus \text{\circled{4}}
 \end{pmatrix}, \label{cpo1}
\end{gather}
where
\begin{gather*}
 \text{\circled{1}} = \phi_1 v_2 - \phi_2 v_1 + a_2w_1 -a_1w_2, \\
 \text{\circled{2}}= \langle\!\langle w_2,v_1\rangle\!\rangle - \langle\!\langle w_1,v_2\rangle\!\rangle,\\
 \text{\circled{3}}= [\phi_1,\phi_2]-w_2^\flat \otimes w_1 + w_1^\flat \otimes w_2 - sv_2^\flat \otimes v_1 + sv_1^\flat \otimes v_2, \\
 \text{\circled{4}}= \phi_1 w_2 - \phi_2 w_1 + s a_1v_2 -s a_2v_1.
\end{gather*}
Recall that $s=1,0$ or $-1$ according to whether $M=S^n$, ${\mathbb R}^n$ or ${\mathbb H}^n$ respectively. The canonical flat connection $\nabla$ on $E$ is given by
\begin{gather}
 \nabla_u
 \begin{pmatrix}
 v \oplus a \\ \phi \oplus w
 \end{pmatrix}=
 \begin{pmatrix}
 \text{\circled{1}}\oplus\text{\circled{2}} \\ \text{\circled{3}}
 \oplus \text{\circled{4}}
 \end{pmatrix}, \label{cpo2}
\end{gather}
where
\begin{gather*}
 \text{\circled{1}} = \nabla^\Sigma_uv + \phi u - a \ii(u)^\sharp,\\
 \text{\circled{2}} = {\rm d}a(u) + \ii(u,v) -\langle\!\langle u,w\rangle\!\rangle,\\
 \text{\circled{3}} = \nabla^\Sigma_u \phi + w^\flat \otimes \ii(u)^\sharp - \ii(u)\otimes w +s\big(v^\flat \otimes u - u^\flat \otimes v\big),\\
 \text{\circled{4}}= \nabla^\Sigma w + \phi\big(\ii(u)^\sharp\big) + s au.
\end{gather*}
We may now also rewrite formula \eqref{seven}, expressing the fundamental representation of
\begin{gather}
A(f) \cong
\begin{matrix}
 T \Sigma \\ \oplus \\ \mathfrak{so}(T \Sigma) \oplus T \Sigma
\end{matrix}\label{af}
\end{gather}
on $E$, as
\begin{gather*}
 \bar\nabla_{\left(
 \begin{smallmatrix}
 u \\ \phi \oplus w
 \end{smallmatrix}\right)}
 \begin{pmatrix}
 v \oplus a \\ \psi \oplus x
 \end{pmatrix}=
 \begin{pmatrix}
 \text{\circled{1}}\oplus\text{\circled{2}} \\ \text{\circled{3}}
 \oplus \text{\circled{4}}
 \end{pmatrix},
\end{gather*}
where
\begin{gather*}
 \text{\circled{1}} =\nabla_u^\Sigma v + \phi v + a\big(w - \ii(u)^\sharp\big),\\
 \text{\circled{2}} ={\rm d}a(u) - \langle\!\langle w-\ii(u)^\sharp,v\rangle\!\rangle, \\
 \text{\circled{3}}= \nabla^\Sigma_u \psi + [\phi,\psi]+ \big(w-\ii(u)^\sharp\big)^\flat \otimes x - x^\flat \otimes \big(w - \ii(u)^\sharp\big),\\
 \text{\circled{4}}= \nabla_u^\Sigma x + \phi x - \psi\big(w - \ii(u)^\sharp\big).
\end{gather*}
Under the identification in \eqref{af}, formula \eqref{logs:1} for the Lie bracket on $A(f)$ can be written
\begin{gather}
 \begin{bmatrix}
 \begin{pmatrix}
 v_1 \\ \phi_1 \oplus w_1
 \end{pmatrix},
 \begin{pmatrix}
 v_2 \\ \phi_2 \oplus w_2
 \end{pmatrix}
 \end{bmatrix}=
 \begin{pmatrix}
 [v_1,v_2] \\ \text{\circled{1}} \oplus \text{\circled{2}}
 \end{pmatrix},\label{forbracket}
\end{gather}
where
\begin{gather*}
 \text{\circled{1}}= [\phi_1,\phi_2] + \nabla^\Sigma_{v_1}\phi_2 - \nabla^\Sigma_{v_2}\phi_1 + w_1^\flat \otimes w_2 - w_2^\flat \otimes w_1 \\
\hphantom{\text{\circled{1}}=}{} + w_2^\flat \otimes \ii(v_1)^\sharp - \ii(v_1) \otimes w_2 - w_1^\flat \otimes \ii(v_2)^\sharp + \ii(v_2) \otimes w_1 + s\big(v_2^\flat \otimes v_1 - v_1^\flat \otimes v_2\big),\\
 \text{\circled{2}}= \nabla^\Sigma_{v_1}w_2 + \phi_1 w_2 - \nabla^\Sigma_{v_2}w_1 -\phi_2 w_1 - \phi_1\big(\ii(v_2)^\sharp\big) + \phi_2\big(\ii(v_1)^\sharp\big).
\end{gather*}

\subsection*{The Gauss--Codazzi equations}
While the connection $\nabla $ on $E$ is flat (it corresponds to the canonical flat connection on ${\mathfrak g} \times \Sigma $), this fact is hidden in the present model of $E$. So we compute its curvature:
\begin{gather}
 \curvature \nabla(u_1,u_2)
 \begin{pmatrix}
 v \oplus a \\ \phi \oplus w
 \end{pmatrix}=
 \begin{pmatrix}
 \text{\circled{1}}\oplus\text{\circled{2}} \\ \text{\circled{3}}
 \oplus \text{\circled{4}}
 \end{pmatrix},\label{glod}
\end{gather}
where
\begin{gather*}
 \text{\circled{1}}= \Xi(v) - a \Theta^\sharp,\\
 \text{\circled{2}}= \Theta(v),\\
 \text{\circled{3}}= \adjoint_\phi \Xi + w^\flat \otimes\Theta^\sharp - \Theta \otimes w,\\
 \text{\circled{4}}= \Xi(w) + \phi\big(\Theta^\sharp\big)
\end{gather*}
and
\begin{gather*}
 \Xi = \curvature \nabla^\Sigma (u_1,u_2) +\ii(u_1)\otimes\ii(u_2)^\sharp - \ii(u_2)\otimes\ii(u_1)^\sharp + s\big(u_1^\flat \otimes u_2 - u_2^\flat \otimes u_1\big),\\
 \Theta =\big(\nabla_{u_1}^\Sigma\ii\big)u_2 - \big(\nabla^\Sigma_{u_2}\ii\big)u_1.
\end{gather*}
Here $(\adjoint_\phi \Xi) (u):=\phi (\Xi (u)) - \Xi(\phi u)$. The equations $\Xi=0$ and $\Theta =0$ are the well-known {\df Gauss--Codazzi equations}. Putting $\curvature \nabla =0$ we obtain the following well-known fact:
\begin{Proposition}\label{gcsecP}
 For any immersed submanifold $\Sigma \subset M$ the Gauss--Codazzi equations hold:
 \begin{gather}
 \curvature \nabla^\Sigma (u_1,u_2) +\ii(u_1)\otimes\ii(u_2)^\sharp - \ii(u_2)\otimes\ii(u_1)^\sharp + s\big(u_1^\flat \otimes u_2 - u_2^\flat \otimes u_1\big)=0,\nonumber\\
 \big(\nabla_{u_1}^\Sigma\ii\big)u_2 - \big(\nabla^\Sigma_{u_2}\ii\big)u_1 =0.\label{gc}
 \end{gather}
\end{Proposition}

\subsection*{The invariant filtration}
We now turn attention to the first three terms
\begin{gather}
 A(f)\supset A(f)^2 \supset A(f)^3\label{jjtt}
\end{gather}
of the invariant filtration associated with the immersion $f \colon \Sigma \rightarrow M$. Let $\mathfrak{so}[T \Sigma]_\ii$ denote the isotropy of $\ii $ under the representation of $\mathfrak{so}[T \Sigma]$ on~$S^2\,T^*\Sigma$ induced by the canonical representation on $T \Sigma $, so that we have a filtration
\begin{gather}
 \mathfrak{so}[T^+\Sigma] \supset \mathfrak{so}[T \Sigma] \supset \mathfrak{so}[T \Sigma]_\ii \label{modelfil}.
\end{gather}
Note that $\mathfrak{so}[T^+\Sigma]_\ii$ is generally intransitive.
\begin{Proposition}\label{ifP} Define a monomorphism ${i} \colon \mathfrak{so}[T^+\Sigma] \rightarrow E$ by
 \begin{gather}
 \begin{pmatrix}
 u \\ \phi
 \end{pmatrix}
 \oplus w \mapsto
 \begin{pmatrix}
 u \oplus 0\\ \phi \oplus (\ii(u)^\sharp + w)
 \end{pmatrix}.\label{mon}
 \end{gather}
 Then:
 \begin{enumerate}\itemsep=0pt
 \item[$1.$]
 The monomorphism is a Lie algebroid isomorphism of $\mathfrak{so}[T^+\Sigma]$ onto $A(f)$ that maps the filtration~\eqref{modelfil} to the filtration~\eqref{jjtt}.
 \item[$2.$] 
 Under the corresponding identification $A(f) \!\cong\! \mathfrak{so}[T^+\Sigma]$ the canonical representation of~$A(f)$ on $T_\Sigma M \cong T^+ \Sigma $ described in Remark~{\rm \ref{cR}} coincides with the canonical representation of~$\mathfrak{so}[T^+\Sigma] $ on~$T^+ \Sigma $.
 \item[$3.$] 
 Under the corresponding identification $A^2(f) \cong \mathfrak{so}[T\Sigma]$ the canonical representation \linebreak of~$A^2(f)$ on~$T \Sigma $ described in Proposition~{\rm \ref{cP}} coincides with the canonical representation of~$\mathfrak{so}[T\Sigma]$ on~$T \Sigma $.
 \end{enumerate}
\end{Proposition}

\begin{Remark}\label{ifR} Note that by (3), and the fact that the metric on $T \Sigma $ is the {\em unique} $\mathfrak{so}[T \Sigma]$-invariant metric up to scale, we recover the metric on $\Sigma$ inherited from $M$, up to scale, as part of the intrinsic geometry on $\Sigma $ defined by the purely abstract considerations of Section~\ref{invint}. (It is also, however, encoded in `extrinsic' data! See Propositions~\ref{ppmP} and~\ref{pzeroP} below.)
\end{Remark}

\begin{proof}[Proof of proposition] It is easy to see that ${i}$ has $A(f)$ as its image. To show that ${i} $ is a Lie algebroid morphism is a straightforward computation using the formula~\eqref{forbracket} for the bracket
 on $A(f)$ and the Gauss--Codazzi equations~\eqref{gc}.

 We now compute $A^2(f)$ and $A^3(f)$ using the isotropy characterization of Proposition~\ref{cP}. First note that the kernel of the anchor of $A(f)$ is
 \begin{gather*}
 {\mathfrak h} = \begin{matrix}
 0 \\ \oplus \\ \mathfrak{so}(T \Sigma) \oplus T \Sigma,
 \end{matrix}
 \end{gather*}
 so that $E/{\mathfrak h} \cong T \Sigma \oplus ({\mathbb R} \times \Sigma)$. Under this identification, the natural inclusion $T \Sigma \cong A(f)/{\mathfrak h} \hookrightarrow E/{\mathfrak h} $ is given by $v \mapsto v \oplus 0$. The representation of $A(f)$ on $E/{\mathfrak h}$ induced by the fundamental representation is readily calculated and given by
\begin{gather*}
 \bar\nabla_{\left(
 \begin{smallmatrix}
 u\\ \phi \oplus w
 \end{smallmatrix}\right)}
 \begin{pmatrix}
 v \oplus a
 \end{pmatrix}
 =\begin{pmatrix}
 \text{\circled{1}} \oplus \text{\circled{2}}
 \end{pmatrix},
\end{gather*}
where
\begin{gather*}
 \text{\circled{1}}= \nabla^\Sigma_u v + \phi v + a\big(w -\ii(u)^\sharp\big) \qquad\text{and}\qquad \text{\circled{2}} ={\rm d}a(u) - \langle\!\langle w - \ii(u)^\sharp,v\rangle\!\rangle.
\end{gather*}
It is not too hard to see that the isotropy $A^2(f) \subset A(f)$ of $T\Sigma \subset E/{\mathfrak h}$ under this representation is ${i}(\mathfrak{so}[T\Sigma])$. A routine computation now also establishes~(2) and~(3).

The kernel of the anchor of $A^2(f) \subset A(f)$ is
\begin{gather*}
 {\mathfrak h}^2 = \begin{matrix}
 0 \\ \oplus \\ \mathfrak{so}(T \Sigma) \oplus 0,
 \end{matrix}
\end{gather*}
so that
\begin{gather*}
 A(f)/{\mathfrak h}^2 \cong
 \begin{matrix}
 T \Sigma \\ \oplus \\ T \Sigma.
 \end{matrix}
\end{gather*}
The canonical inclusion $T \Sigma \hookrightarrow A(f)/{\mathfrak h}^2$ of Proposition~\ref{cP} is given by
 \begin{gather}
 v \mapsto \begin{pmatrix}
 v \\ \ii(v)^\sharp
 \end{pmatrix}.\label{vue}
 \end{gather}
 The representation of $A^2(f) \cong \mathfrak{so}[T\Sigma]=T \Sigma \oplus \mathfrak{so}(T \Sigma)$ on $A(f)/{\mathfrak h}^2$, induced by the fundamental representation, is now determined to be
\begin{gather*}
\bar\nabla_{\left(
 \begin{smallmatrix}
 u \\ \phi
 \end{smallmatrix}\right)}
\begin{pmatrix}
 v \\ x
\end{pmatrix}=
\begin{pmatrix}
 \nabla^\Sigma_u v + \phi v \\ \nabla^\Sigma_u x + \phi x
\end{pmatrix}.
\end{gather*}
That is, it is just two copies of the canonical representation of~$\mathfrak{so}[T\Sigma]$ on $T \Sigma $. One computes
\begin{gather*}
 \bar\nabla_{\left(
 \begin{smallmatrix}
 u \\ \phi
 \end{smallmatrix}\right)}
\begin{pmatrix}
 v \\ \ii(v)^\sharp
\end{pmatrix}=
\begin{pmatrix}
 \nabla^\Sigma_u v + \phi v \\ \ii\big(\nabla^\Sigma_u v + \phi v\big)^\sharp
\end{pmatrix}+
\begin{pmatrix}
 0 \\ \big( \big(\nabla^\Sigma_u \ii + \adjoint_\phi
 \ii\big)(v) \big)^\sharp
\end{pmatrix},
\end{gather*}
where $(\adjoint_\phi\ii)(v_1,v_2):=-\ii(\phi v_1,v_2)-\ii(v_1,\phi v_2)$. But the representation of ${\mathfrak{so}[T\Sigma]} $ on $S^2\,T^* \Sigma$ determined by the canonical representation is given by
\begin{gather*}
 \bar\nabla_{\left(
 \begin{smallmatrix}
 u \\ \phi
 \end{smallmatrix}\right)} \ii=\nabla^\Sigma_u \ii + \adjoint_\phi \ii.
\end{gather*}
Keeping the formula \eqref{vue} for the inclusion $T \Sigma \hookrightarrow A(f)/{\mathfrak h}^2$ in mind, we see that~$A(f)^3$~-- the isotropy of $T\Sigma \subset A(f)^2/{\mathfrak h}$ under the above representation of
$A(f)^2 \cong \mathfrak{so}[T\Sigma] $ -- is $\mathfrak{so}[T\Sigma]_\ii$, which completes the proof of~(1).
\end{proof}

\section{The Bonnet theorem and its proof}\label{bonnetsec}
We continue to let $M$ denote $S^n$, ${\mathbb R}^n$, or ${\mathbb H}^n $ (with $s=1,0,-1$) as in the preceding section. Our first task is to explain how the first and second fundamental forms of an immersion into $M$ can be recovered from more abstract considerations, preparing the way for the abstract formulation of the Bonnet theorem. The classical analogue of this result is then proven as a corollary.

\subsection*{A polynomial invariant for elliptic and hyperbolic subgeometry}
We now fix an appropriate multiple of the Killing form on~${\mathfrak g} $, relevant to the cases $M=S^n$ and ${\mathbb H}^n$. To this end, suppose ${B} $ is a vector bundle equipped with an inner product $\langle\!\langle \,\cdot\,,\,\cdot\,\rangle\!\rangle$ and let $U \mapsto U^\flat$ denote the corresponding isomorphism $ {B} \rightarrow {B}^*$. On the Lie algebra bundle $\mathfrak{so}({B} )$ of skew-symmetric endomorphisms of fibres of~${B} $ is a natural inner product $q^{B} $. If we write $U\wedge V := U^\flat \otimes V - V^\flat \otimes U$, then $q^{B} $ is given by
\[q^{B} (U_1\wedge V_1, U_2\wedge V_2)= \langle\!\langle U_1,U_2\rangle\!\rangle \langle\!\langle V_1,V_2\rangle\!\rangle-\langle\!\langle U_1,V_2\rangle\!\rangle \langle\!\langle U_2,V_1\rangle\!\rangle. \]
This form is a positive definite multiple of the Killing form, the multiplicative factor depending only on $\rank {B}$. We record the identity
\begin{gather}
 q^{B}(U\wedge V, \Phi)=\langle\!\langle\Phi\, U,V\rangle\!\rangle,\label{vege}
\end{gather}
holding for any sections $U$, $V$ of $ {B} $ and $\Phi$ of $\mathfrak{so}({B})$.
\begin{Lemma}\label{ppmL} For $s= \pm 1$ there exists a unique multiple $Q$ of the Killing form on ${\mathfrak g} $ with the following property: If the corresponding bilinear form on the trivial bundle ${\mathfrak g} \times M \cong \mathfrak{so}[TM]=TM \oplus \mathfrak{so}(TM)$ is also denoted $Q$, then
 \begin{gather}
 Q\left(
 \begin{pmatrix}
 V_1 \\ \Phi_1
 \end{pmatrix},
 \begin{pmatrix}
 V_2 \\ \Phi_2
 \end{pmatrix}
 \right)=\langle\!\langle V_1,V_2\rangle\!\rangle + s q^{T M}(\Phi_1,\Phi_2), \label{petshop}
 \end{gather}
 where $\langle\!\langle \,\cdot\,,\,\cdot\,\rangle\!\rangle$ is the metric on $M$. Moreover, for any $m \in M$ the isotropy algebra ${\mathfrak g}_m$ intersects its $Q$-orthogonal complement ${\mathfrak g}_m^\perp$ transversally, with the restriction of $Q$ to ${\mathfrak g}_m^\perp$ being positive definite, so that $Q$ drops to a positive definite quadratic form on ${\mathfrak g} /{\mathfrak g}_m$.
\end{Lemma}
\begin{proof} Recall that under the isomorphism ${\mathfrak g}\times M \cong \mathfrak{so}[TM]$ the canonical flat connection on ${\mathfrak g} \times M$ is represented by the Cartan connection $\nabla $ given by equations~\eqref{form} and \eqref{torsion2}, and that the canonical algebraic bracket on ${\mathfrak g} \times M$ is represented by the bracket on $\mathfrak{so}[TM]$ given by \eqref{torsion} and \eqref{torsion2}. Since the quadratic form $Q$ defined by~\eqref{petshop} is seen to be $\nabla$-parallel, it follows that it is defined by a quadratic form on ${\mathfrak g} $ (because $M$ is simply-connected). Since $Q$ is $\adjoint$-invariant on each fibre of $\mathfrak{so}[TM]$ with respect to the Lie bracket defined by the algebraic bracket on $\mathfrak{so}[TM]$ (verify this with the help of~\eqref{vege}) and is also non-degenerate, the corresponding quadratic form on $\mathfrak g$ must be a multiple of the Killing form, because ${\mathfrak g} $ is semisimple.

Fixing $m \in M$, the identification ${\mathfrak g} \times M \cong \mathfrak{so}[T\Sigma]=TM \oplus \mathfrak{so}(TM)$ gives us an identification of ${\mathfrak g} $ with $T_mM \oplus \mathfrak{so}(T_mM)$, in which the isotropy ${\mathfrak g}_m$ is identified with the second summand. The remaining claims of the lemma now follow from~\eqref{petshop}.
\end{proof}

We will use $Q$ to denote both the multiple of the Killing form on $\mathfrak g$ and the corresponding form given by~\eqref{petshop}. Continuing to suppose $s=\pm 1$, let $\omega \colon A \rightarrow {\mathfrak g} $ be an infinitesimal immersion of an $(n-1)$-dimensional manifold $\Sigma$ into $M$. The above polynomial invariant $Q$ defines a~quadratic form on $E={\mathfrak g} \times \Sigma$ that we also denote by~$Q$ in the following.

According to Axiom I3, and its consequence \eqref{worm}, when we regard the kernel ${\mathfrak h} $ of the anchor $\# \colon A \rightarrow T \Sigma $ as a subbundle of $E = {\mathfrak g} \times \Sigma $ using $\omega $, the fibre of ${\mathfrak h}|_x $ over $x \in \Sigma $ becomes identified with ${\mathfrak g}_m \times \{x\}$, for some $m \in M$. It follows from the second part of the lemma above that the $Q$-orthogonal complement ${\mathfrak h}^\perp $ intersects ${\mathfrak h} $ transversally and that $Q$ drops to a positive definite quadratic form on~$E/{\mathfrak h} $. The restriction of this form to $T \Sigma \subset E/{\mathfrak h} $ is a metric on $\Sigma $ we denote by $\langle\!\langle \,\cdot\,,\,\cdot\,\rangle\!\rangle_\omega$ and call the {\df first fundamental form} of $\omega$. This form is evidently an invariant of $\omega $.

\begin{Proposition}\label{ppmP} If $s=\pm 1$ and $f \colon \Sigma \rightarrow M$ is an infinitesimal immersion, then $\langle\!\langle\,\cdot\,,\,\cdot\,\rangle\!\rangle_{\delta f}$ coincides with the inherited metric.
\end{Proposition}
\begin{proof} Using \eqref{petshop} we can derive a formula for the quadratic form $Q$ on $E$ with respect to the model of $E$ given in~\eqref{mmmm}:
\begin{gather}
 Q\left(
 \begin{pmatrix}
 v_1 \oplus a_1 \\ \phi_1 \oplus w_1
 \end{pmatrix},
 \begin{pmatrix}
 v_2 \oplus a_2 \\ \phi_2 \oplus w_2
 \end{pmatrix}
 \right) = \langle\!\langle v_1,v_2\rangle\!\rangle + a_1 a_2
 + s q^{T\Sigma}(\phi_1,\phi_2) + s\langle\!\langle
 w_1,w_2\rangle\!\rangle,\label{foot}
\end{gather}
where $\langle\!\langle \,\cdot \,,\,\cdot\,,\rangle\!\rangle$ is the metric on $\Sigma$. The proposition follows.
\end{proof}

\subsection*{A polynomial invariant for Euclidean subgeometry}
In the case $s=0$ the Killing form is less useful and we instead fix an $\Adjoint$-invariant quadratic form on the radical\footnote{In the special case $n=2$, we use the commutator $[{\mathfrak g},{\mathfrak g} ]$ instead of $\rad {\mathfrak g} $.} $\rad {\mathfrak g}$ (a form which is {\em not} however invariant under outer automorphisms of $\mathfrak g$). Note that under the identification ${\mathfrak g} \times M \cong TM \oplus \mathfrak{so}(TM)$ delivered by Proposition \ref{ambientP}, $\rad {\mathfrak g} \times M $ coincides with $TM \oplus 0 \cong TM$. Viewing elements of $\rad{\mathfrak g} $ concretely as constant vector fields on $M$, the isomorphism $\rad {\mathfrak g} \times M \rightarrow TM $ is given simply by $(\xi,m)\mapsto \xi(m)$.
\begin{Lemma}\label{pzeroL} If $s=0$ then there exists a unique $\Adjoint$-invariant positive definite quadratic form~$Q$ on $\rad {\mathfrak g}$ such that the corresponding quadratic form on $\rad {\mathfrak g} \times M \cong TM$ is the metric on~$M$. Furthermore, for any $m \in M$, the isotropy subalgebra ${\mathfrak g}_m$ has the following property: the restriction of the projection ${\mathfrak g} \rightarrow {\mathfrak g} /{\mathfrak g}_m$ to $\rad {\mathfrak g} $ is an isomorphism.
\end{Lemma}

The easy proof of this lemma is omitted. Continuing to suppose $s=0$, let $\omega \colon A \rightarrow {\mathfrak g} $ be an infinitesimal immersion of an $n$-dimensional manifold $\Sigma$ into $M$. The above quadratic form $Q$ on $[{\mathfrak g}, {\mathfrak g}] $ defines an positive definite quadratic form on $\rad E:=\rad {\mathfrak g} \times \Sigma$, also denoted here by~$Q$.

Once again, let $\omega \colon A \rightarrow {\mathfrak g} $ be an infinitesimal immersion of an $(n-1)$-dimensional manifold~$\Sigma$ into $M$, with $s=0$, and let ${\mathfrak h} $ be the kernel of the anchor $\# \colon A \rightarrow T \Sigma $. According to Axiom~I3, when we regard ${\mathfrak h} $ as a subbundle of $E = {\mathfrak g} \times \Sigma $ using $\omega $, the fibre of ${\mathfrak h}$ over $x \in \Sigma $ becomes identified with
${\mathfrak g}_m \times \{x\}$, for some $m \in M$. It follows from the second part of the lemma above that the projection $E \rightarrow E/{\mathfrak h} $ restricts to an isomorphism $\rad E \rightarrow E/{\mathfrak h} $ pushing the quadratic form~$Q$ on $\rad E$ forward to a positive definite quadratic form on~$E/{\mathfrak h} $. The restriction of this form to $T \Sigma \subset E/{\mathfrak h} $ is a metric on $\Sigma $ we denote by $\langle\!\langle \,\cdot\,,\,\cdot\,\rangle\!\rangle_\omega$ and call the {\df first fundamental form} of $\omega$. This form is evidently an invariant of~$\omega $.

\begin{Proposition}\label{pzeroP} If $s=0$ and $f \colon \Sigma \rightarrow M$ is an infinitesimal immersion, then $\langle\!\langle\,\cdot\,,\,\cdot\,\rangle\!\rangle_{\delta f}$ coincides with the inherited metric on~$\Sigma $.
\end{Proposition}
\begin{proof} In terms of the model of $E$ given in \eqref{mmmm}, we have
 \begin{gather*} \rad E= \begin{matrix}
 T \Sigma \oplus 0 \\ \oplus \\ 0 \oplus 0,
 \end{matrix}
 \end{gather*}
 and the inner product on $\rad E$ defined by $Q$ is given by
 \begin{gather*}
 Q\left( \begin{pmatrix}
 v_1 \oplus 0\\ 0\oplus 0
 \end{pmatrix},
 \begin{pmatrix}
 v_2 \oplus 0 \\ 0 \oplus 0
 \end{pmatrix} \right)= \langle\!\langle v_1,v_2\rangle\!\rangle,
 \end{gather*}
 where $\langle\!\langle \,\cdot \,,\,\cdot\,,\rangle\!\rangle$ is the metric on $\Sigma$. The proposition easily follows.
\end{proof}

\subsection*{The second fundamental form of an infinitesimal immersion}
Let $\omega \colon A \rightarrow {\mathfrak g} $ be an infinitesimal immersion of an $(n-1)$-dimensional manifold $\Sigma$ into $M$, with $s \in \{-1,0,1\}$. Assuming for the moment that $M $ is not an even-dimensional sphere, Proposition~\ref{orientedP} may be applied to orient $E/{\mathfrak h} $ in a natural way.

The {\df normal} $\xi $ of $\omega $ is a section of $E={\mathfrak g} \times \Sigma$ defined as follows. In the case $s=\pm 1$, the $Q$-orthogonal complement $A^\perp$ of $A$ in $E$ has rank one and $A^\perp \subset {\mathfrak h}^\perp$, where ${\mathfrak h} $ is the kernel of the anchor $\# \colon A \rightarrow T \Sigma $. In particular, since the restriction of $Q$ to ${\mathfrak h}^\perp$ is positive definite, we must have $A^\perp \cap A =0$ and non-zero elements of $A^\perp$ must have positive $Q$-length. The image of such an element under the projection $E \rightarrow E/{\mathfrak h} $ is necessarily transverse to $T \Sigma $ (the image of $A$ under the projection). It follows that there exists a unique section $\xi$ of $A^\perp$ with constant $Q$-length one, whose image under the natural projection $E \rightarrow E/{\mathfrak h}$ is transverse to $T \Sigma $ and has sense consistent with the orientations of $T \Sigma $ and~$E/{\mathfrak h}$.

Recall that in the $s=0$ case, the natural projection $E \rightarrow E/{\mathfrak h} $ restricts to an isomorphism $\rad E \rightarrow E/{\mathfrak h} $. Since $T \Sigma$ is the image of $A$ under the projection $E \rightarrow E/{\mathfrak h} $, it follows that $A \cap \rad E$ has corank one in $\rad E$ because $T \Sigma $ has corank one in $E/{\mathfrak h} $. Therefore there exists a~unique section~$\xi $ of~$\rad E$ $Q$-orthogonal to $A \cap \rad E$ that has constant $Q$-length one, and whose image under the natural projection $E \rightarrow E/{\mathfrak h}$ is transverse to $T \Sigma $ and has sense consistent with the with the orientations of~$T \Sigma $ and $E/{\mathfrak h}$.

For any $s \in \{-1,0,1\}$ we have
\begin{Lemma}\label{iiL} $\nabla_u\xi \in A $ for all $u \in T \Sigma $.
\end{Lemma}
\begin{proof} Differentiate the requirement $Q(\xi,\xi)=1$ and appeal to the fact that $Q$ is $\nabla$-parallel and, in the case $s=0$, that $\rad E$ is $\nabla $-invariant.
\end{proof}

The {\df second fundamental form} $\ii_\omega$ of $\omega $ is the two-tensor on $\Sigma $ defined by
\begin{gather*}
 \ii_\omega(u,v)=-\langle\!\langle u,\# \nabla _v \xi \rangle\!\rangle_\omega.
\end{gather*}
\begin{Proposition}\label{iiP} The second fundamental form of the logarithmic derivative of an isometric immersion $f \colon \Sigma \rightarrow M$ coincides with the second fundamental form of $f$ in the usual sense:
 $\ii_{\delta f}=\ii $.
\end{Proposition}
\begin{proof} With respect to the model of $E$ in \eqref{mmmm}, the normal for the logarithmic derivative is given by
 \begin{gather}
 \xi = \begin{pmatrix}
 0 \oplus 1 \\ 0 \oplus 0 \label{n}
 \end{pmatrix}.
 \end{gather}
 The proof is now a straightforward computation.
\end{proof}

If $M$ is an even-dimensional sphere, then $E/{\mathfrak h} $ is still orientable if we assume $\Sigma $ is simply-connected, but there is no canonical orientation. Fixing the orientation arbitrarily, the second fundamental form obtained is an invariant up to sign only. However, we must remember that this anomaly is explained by the fact that the symmetries in this case include the orientation-reversing isometries.

\subsection*{The abstract Bonnet theorem}
Recall that $M = S^n$, ${\mathbb R}^n$ or ${\mathbb H}^n $, and that $s \in \{1,0,-1\}$ is the corresponding value of the scalar curvature. We are viewing $M$ as a Klein geometry with transitively acting group $G$ the group of orientation-preserving isometries, and identify the Lie algebra ${\mathfrak g} $ of $G$ with the Killing fields on $M$. For the moment $\Sigma $ denotes a smooth $(n-1)$-dimensional manifold without additional structure.
\begin{Theorem}\label{atlastT} Let $\omega \colon A \rightarrow {\mathfrak g} $ be any infinitesimal immersion of $\Sigma $ into $M$. Then, assu\-ming~$\Sigma$ is simply-connected: {\normalfont (i)} The second fundamental form $\ii_\omega $ is symmetric; and {\normalfont (ii)} $\omega$ admits an immersion $f \colon \Sigma \rightarrow M$ as primitive, whose first and second fundamental forms are precisely $\langle\!\langle \,\cdot\,,\,\cdot\,,\rangle\!\rangle_\omega $ and $\ii_\omega $. The immersion is unique up to a symmetry, i.e., up to an orientation-preserving isometry, unless~$M$ is an even-dimensional sphere, in which case uniqueness is only up to isometry.
\end{Theorem}
\begin{Remark}\label{atlastR} One can easily argue for uniquess up to an orientation-preserving isometry in the case of an even-dimensional sphere under extra hypotheses on $\ii_\omega $ -- for example, if $\ii_\omega $ is sign-definite.
\end{Remark}
\begin{proof} Clearly (i) follows from (ii). Applying Theorem \ref{infirmT}, we obtain an immersion $f \colon M \rightarrow \Sigma $ as primitive for $\omega $, unique up to orientation-preserving isometries (arbitrary isometries in the case of even-dimensional spheres). By definition, this means $\delta f $ and $\omega $ are isomorphic in $\mathbf{Inf}$, and so have the same invariants. This, in turn, implies, by the preceding three propositions, that the first and second fundamental forms of $\omega$ are the first and second fundamental forms of the immersion $f$. In the case $M$ is an even-dimensional sphere we can only guarantee that the second fundamental form of $f$ is $\pm \ii_\omega$. However, a wrong sign is readily corrected for by composing $f$ with the map $m \mapsto -m$.
\end{proof}

\subsection*{The classical Bonnet theorem}
As a corollary we now prove the following very well-known converse of Proposition \ref{gcsecP}:
\begin{Theorem}\label{cbtT} Suppose $\Sigma $ is an $(n-1)$-dimensional, oriented, simply-connected, Riemannian manifold, supporting a section $\ii $ of $S^2\,T^* \Sigma$ satisfying the Gauss--Codazzi equations~\eqref{gc}, with $s \in \{-1,0,1\}$. Then there exists an isometric immersion $f \colon \Sigma \rightarrow M$, where $M$ is respective\-ly~${\mathbb H}^n$,~${\mathbb R}^n $ or~$S^n$, whose second fundamental form is~$\ii$. The immersion is unique up to symmetry $($see the preceding theorem$)$.
\end{Theorem}
\begin{proof} By the preceding theorem it suffices to construct an infinitesimal immersion $\omega \colon A \rightarrow {\mathfrak g} $ whose first fundamental form is the metric on $\Sigma$, and whose second fundamental form is~$\ii $. Applying Construction Principle~\ref{constructionPrinciple}, we define a vector bundle $E$ over~$\Sigma $ by
\begin{gather*}
 E= \begin{matrix}
 T \Sigma \oplus ({\mathbb R} \times \Sigma) \\ \oplus \\
 \mathfrak{so}(T \Sigma) \oplus T \Sigma.
 \end{matrix}
\end{gather*}
This bundle becomes a ${\mathfrak g}$-bundle (${\mathfrak g} $ the Lie algebra of Killing fields of the ambient space $M$) if we define an algebraic bracket on $E$ by~\eqref{cpo1}. The connection $\nabla $ on $E$ defined by~\eqref{cpo2} has curvature given by~\eqref{glod}, by an identical calculation. The Gauss--Codazzi equations ensure that $\nabla $ is flat. The reader will also verify that the algebraic bracket on $E$ is $\nabla$-parallel.

Realize the Lie algebroid $\mathfrak{so}[T^+ \Sigma]$ as a subbundle of $E$ using the monomorphism ${i} \colon \mathfrak{so}[T^+ \Sigma]$ $\rightarrow E$ defined by \eqref{mon}. With the help of the Gauss--Codazzi equations, one shows that the Maurer--Cartan equations~\eqref{mce} hold. The relevant calculations amount to those already used to derive the Lie algebroid bracket on $A(f)$ in the case of an immersion $f$, and to show ${i} $ was a Lie algebroid morphism onto $A(f)$ in that case. Applying Proposition~\ref{croonP}, we obtain a Lie algebroid morphism $\omega \colon \mathfrak{so}[T^+ \Sigma]\rightarrow {\mathfrak g} $. It is not hard to verify that $\omega $ satisfies the axioms of an infinitesimal immersion.

Suppose $s=\pm 1$. Since the quadratic form $Q$ on ${\mathfrak g} $ fixed in Lemma~\ref{ppmL} is an outer invariant of the Lie algebra ${\mathfrak g} $, the corresponding quadratic form on~$E$ must be given by the same formula~\eqref{foot} we derived for immersions, because the formula for the algebraic bracket on $E$ is the same. It follows from its definition that $\langle\!\langle \,\cdot\,,\,\cdot\,\rangle\!\rangle_\omega$ is the prescribed metric on $\Sigma $.

In the $s=0$ case, the quadratic form $Q$ on $\rad E$ is only an inner invariant and the same argument won't work. However, we know that $\langle\!\langle \,\cdot\,,\,\cdot\,\rangle\!\rangle_\omega$ is $A^2$-invariant and, by the same arguments used for immersions, we have $A^2 \cong \mathfrak{so}[T \Sigma] $ acting tautologically on $T \Sigma $. Since all metrics on $T \Sigma $ invariant under this tautological representation coincide with the metric on $T \Sigma $, up to scale, we conclude that $\langle\!\langle \,\cdot\,,\,\cdot\,\rangle\!\rangle_\omega$ is a~positive constant multiple of the metric on $\Sigma $. Now let $\lambda > 0$ be arbitrary. Then it is an elementary property of ${\mathfrak g} $ that there exists an {\em outer} automorphism of ${\mathfrak g} $ pulling $Q$ back to $\lambda Q$. It follows that by replacing $\omega \colon \mathfrak{so}[T^+\Sigma] \rightarrow {\mathfrak g} $ with its composition with an appropriate outer automorphism of ${\mathfrak g} $, we can arrange that $\langle\!\langle \,\cdot\,,\,\cdot\,\rangle\!\rangle_\omega$ coincides with the metric on $\Sigma $.

By construction, the normal $\xi $ of $\omega $ is given by the same formula \eqref{n} as for the logarithmic derivative of immersions, and by the same calculation as before, we have $\ii_\omega=\ii$.
\end{proof}

\section{Moving frames and curves in the equi-affine plane}\label{appendixA}
From the present point of view, a moving frame is just a convenient device for making computations. Moving frames are particularly well-suited to the study of immersed curves.

\subsection*{Moving frames}
Let $M$ be a Klein geometry with transitively acting Lie group $G$ and $f \colon \Sigma \rightarrow M $ an immersion, so that $A(f)$ may be identified with a subbundle of $E:={\mathfrak g} \times \Sigma $. By a~{\df moving frame} on $\Sigma $ let us mean a smooth map $\tilde f \colon \Sigma \rightarrow G$ such that $f(x)=\tilde f(x) \cdot f(x_0)$, for some $x_0 \in \Sigma $. A moving frame determines a new trivialization of $E$: For each $\xi \in {\mathfrak g} $ we define a corresponding {\df moving section}~$\tilde \xi $ of~$E$ by
 \[ \tilde \xi(x) = \Adjoint_{\tilde f(x)} \xi.\]
 Then the map
\[(x,\xi)\mapsto \big(x, \tilde \xi (x)\big)\colon \ E \rightarrow E\]
is an isomorphism which one hopes will `rectify' the Lie algebroid $A(f) \subset E$ for a suitable choice of frame $\tilde f$. Specifically, one seeks to arrange that $A(f)=\spann\{\tilde \xi_1,\ldots,\tilde \xi_k\}$, for some basis $\{\xi_1,\xi_2,\ldots,\xi_k,\xi_{k+1},\ldots, \xi_N\}$ of $\mathfrak g$. Before illustrating the idea with an example, we record the following key observation:
\begin{Proposition}\label{movingP} Let $\nabla $ be the canonical flat connection on $E = {\mathfrak g} \times M$, and $\{\,\cdot\,,\,\cdot\,\}$ the canonical algebraic bracket determined by the Lie bracket on ${\mathfrak g}$. Then for any choice of moving frame $\tilde f \colon \Sigma \rightarrow G$ we have
 \begin{gather}
 \nabla_U \tilde \xi =\big\{\tilde \xi, \delta \tilde f(U)\big\},\qquad \xi \in {\mathfrak g},\label{curio}
 \end{gather}
 where $\delta \tilde f \colon T \Sigma \rightarrow {\mathfrak g} $ is the ordinary logarithmic derivative of $\tilde f $, i.e., $\delta \tilde f = \tilde f^* \omega_G$, where $\omega_G$ denotes the right-invariant Maurer--Cartan form on $G$.
\end{Proposition}
\begin{proof} This follows from the following general identity for a path $t \mapsto g(t) $ on a Lie group $G$:
 \begin{gather*}
 \frac{{\rm d}}{{\rm d}t}\Adjoint_{g(t)}\xi = \left[\Adjoint_{g(t)}\xi, \omega_G\left(\frac{{\rm d}}{{\rm d}t}g(t)\right)\right],\qquad \xi \in {\mathfrak g}.\tag*{\qed}
 \end{gather*}\renewcommand{\qed}{}
\end{proof}

\subsection*{Curves in the equi-affine plane}
Consider the geometry determined by the group $G$ of all area-preserving affine transformations of the plane $M = {\mathbb R}^2 $. We can recover a well-known characterization of curves in this geomet\-ry~-- see the theorem below~-- by combining the general theory of the present article with a choice of moving frame.

By \cite[Example~2.3(2)]{Blaom_F}, we have $\automorphism(M)=G$. Up to a choice of scale and starting point~$a$, a~suitably regular curve $\Gamma $ has a canonical parameterization $f \colon [a,b] \rightarrow {\mathbb R}^2 $ for which the derivative~$\dot f$ is a parallel vector field on the curve, in the intrinsic geometry $\Gamma \subset {\mathbb R}^2 $ inherits as a~submanifold (see Section~\ref{invint}), this geometry being an infinitesimal
parallelism. Since it is easy to anticipate the form of this parameterization, the preceding statement will be justified ex post facto.

Regard the time derivatives $\dot f(t), \ddot f(t) \in {\mathbb R}^2 $ in some arbitrary parameterization $f$ of $\Gamma $ as column vectors and define a $2 \times 2$ matrix~$R(t)$ by
\begin{gather*}
 R(t) = \begin{pmatrix}
 \dot f(t) & \ddot f(t)
 \end{pmatrix}.
\end{gather*}
We suppose $\Gamma $ is {\df regular} in the sense that $\det R(t)\ne 0$ for some such $f$. Since $G$ preserves areas, it is easy to see that $\sigma(t)=\det R(t)$, known as the {\df equi-affine speed} of the path~$f$, is an invariant of $f$. It is given by
\[\sigma(t)=s(t)^3 \kappa(t),\]
where $s(t)$ and $\kappa(t)$ are the Euclidean speed and curvature of $f$. Supposing $\Gamma $ is regular, we may, by reparameterizing, arrange that $\sigma(t)$ is some constant $\sigma>0 $, which fixes $f \colon I \rightarrow {\mathbb R}^2 $~-- and, in particular, an orientation for $\Gamma$~-- up to a choice of starting point for the interval~$I$.

As we shall prove, an invariant function on $\Gamma $ characterising the curve is the function whose pullback to $I$ by $f$ is given by
\begin{gather*}
 \kk (t) = \frac{1}{\sigma}\det \dot R(t), \qquad \text{where} \quad \dot R(t) =
 \begin{pmatrix}
 \ddot f(t) & \dddot f(t)
 \end{pmatrix}.
\end{gather*}
This function is called the {\df equi-affine curvature} of the curve.

Now the matrices $R(t)$ and $\dot R(t)$ are not independent, for the first column of $\dot R(t) $ is the second column of $R (t)$, implying
 \begin{gather}
 \dot R(t) = R(t) \begin{pmatrix}
 0 & a(t) \\ 1 & b(t)
 \end{pmatrix}, \label{wmm}
 \end{gather}
 for some $a(t)$, $b(t)$. Moreover, one has the following identity for the derivatives of determinants, true for any differentiable path of square matrices $t \mapsto R(t)$:
 \begin{gather*}
 \frac{{\rm d}}{{\rm d}t}\log(\det(R(t))=\trace\big(R(t)^{-1}\dot R(t)\big).
 \end{gather*}
 In the present case the left-hand side vanishes because we assume $\sigma=\det R(t)$ is a constant, and we deduce from~\eqref{wmm} that $b(t)=0$. Taking determinants of both sides of~\eqref{wmm} we obtain $\kk(t)=-a(t)$, and whence arrive at the following formula:
 \begin{gather}
 R(t)^{-1}\dot R(t)= \begin{pmatrix}
 0 & -\kk(t) \\ 1 & 0
 \end{pmatrix}.\label{lemon}
 \end{gather}
 We are now ready to compute the invariant filtration associated with the path $f$.

\begin{Notation}\label{gtttN} For a $2 \times 2$ matrix $A$ we let $A^\dagger$ denote the vector field on ${\mathbb R}^2$ defined by
 \begin{gather*}
 A^\dagger(m)=\frac{{\rm d}}{{\rm d}t} e^{-tA}m\Big|_{t=0}, \qquad m \in {\mathbb R}^2,
 \end{gather*}
 so that $[A^\dagger, B^\dagger]=(AB-BA)^\dagger$. We denote the standard coordinate functions on ${\mathbb R}^2 $ by $\xi$, $\eta$, pulling back to functions $x$, $y $ under the map~$f$, i.e., $f(t)=(x(t), y(t))$.
\end{Notation}

Identifying elements of the Lie algebra ${\mathfrak g} $ of $G$ with vector fields on ${\mathbb R}^2$, we have a basis $\{U,V,X,Y,H\}$ for~${\mathfrak g} $ given by
\begin{gather*}
 U= \frac{\partial}{\partial \xi},\qquad V= \frac{\partial }{\partial \eta},\qquad
 X=\begin{pmatrix} 0 & 1 \\ 0 & 0 \end{pmatrix}^\dagger,\qquad
 Y=\begin{pmatrix} 0 & 0 \\ 1 & 0 \end{pmatrix}^\dagger,\qquad
 H=\begin{pmatrix} 1 & 0 \\ 0 & -1 \end{pmatrix}^\dagger.
\end{gather*}
The bracket on ${\mathfrak g} $ is given by
\begin{alignat*}{6}
& [U,V]=0,\qquad &&[X,U]=0, \qquad && [X,V]=U, \qquad && [Y,U] = V, \qquad && [Y,V]=0,&\\
& [H,U]=U, \qquad && [H,V]=-V, \qquad && [X,Y]=H, \qquad && [H,X]=2X, \qquad && [H,Y]=-2Y.&
\end{alignat*}
A moving frame $\tilde f \colon I \rightarrow G$ is given by
\begin{gather}
 \tilde f (t)(m) = \frac{1}{\sqrt{\sigma }}R(t) m + f(t), \qquad m \in {\mathbb R}^2.\label{formform}
\end{gather}
That is, $\tilde f(t) \in G$ is the affine transformation that sends the origin $(0,0)$ to $f(t)$, and the standard basis of $T_{(0,0)}{\mathbb R}^2$ to the `moving frame' $\big\{\dot f(t)/\sqrt{\sigma}, \ddot f(t)/\sqrt{\sigma}\big\}$ at $T_{f(t)}{\mathbb R}^2 \cong {\mathbb R}^2$. The factor $\sqrt{\sigma}$ ensures that this transformation is area-preserving.

Sections of $ {\mathfrak g}\times I $ are generated by the `moving sections' $\tilde U$, $\tilde V$, $\tilde X$, $\tilde Y$, $\tilde H$, where $\tilde \xi (t) :=\Adjoint_{\tilde f(t)} \xi$. These sections must obey bracket relations analogous to those above:
\begin{gather}
 \{\tilde U,\tilde V\}=0,\qquad \{\tilde X,\tilde U\}=0, \qquad \{\tilde X,\tilde V\}=\tilde U, \qquad \{\tilde Y,\tilde U\} =\tilde V, \qquad \{\tilde Y,\tilde V\}=0,\!\!\!\!\!\! \label{relations}\\
 \{\tilde H,\tilde U\}=\tilde U, \qquad\! \{\tilde H,\tilde V\}=- \tilde V,\qquad\! \{\tilde X,\tilde Y\}= \tilde H, \qquad\! \{\tilde H,\tilde X\}=2\tilde X, \qquad\! \{\tilde
 H,\tilde Y\}=-2\tilde Y. \nonumber
\end{gather}

As we identify an element $\xi $ of ${\mathfrak g} $ with a vector field, $\tilde \xi=\Adjoint_{\tilde f(t)}\xi $ is just the pushforward of $\xi $ under the transformation $\tilde f(t)$. In particular, we see that
\begin{gather*}
 \tilde U(t) = \frac{\dot x(t)}{\sqrt{\sigma}} U + \frac{\dot y(t)}{\sqrt{\sigma}} V, \qquad \tilde V(t) = \frac{\ddot x(t)}{\sqrt{\sigma}} U + \frac{\ddot y(t)}{\sqrt{\sigma}} V,
\end{gather*}
while $\tilde X(t), \tilde Y(t), \tilde H(t) \in {\mathfrak g} $ are all vector fields that vanish at $f(t)$. It follows that $A(f) \subset {\mathfrak g}\times I $ is spanned by the sections $\tilde U$, $\tilde X$, $\tilde Y$, $\tilde H$ (i.e., is `rectified' by our choice of moving frame) and we have
\begin{gather}
 \# \tilde U = \frac{1}{\sqrt{\sigma}}\frac{\partial}{\partial t}, \qquad \#X = \#Y =\# Z =0.\label{anchors}
\end{gather}

The ordinary logarithmic derivative of $\tilde f $ is given by
\begin{gather*}
 \delta \tilde f\left(\frac{\partial}{\partial t}(t)\right)= -\operatorname{Ad}_{\tilde f(t)}\big( R(t)^{-1}\dot
 R(t) \big)^\dagger + \dot x(t) \frac{\partial }{\partial \xi } + \dot y(t) \frac{\partial }{\partial \eta }.
\end{gather*}
Applying the identity \eqref{lemon} above, we obtain
\begin{gather*}
 \delta \tilde f\left(\frac{\partial}{\partial t}(t)\right)= \kk (t) \tilde X(t) - \tilde Y(t) + \sqrt{\sigma}\,\tilde U(t).
\end{gather*}
It follows from \eqref{curio} that
\begin{gather*}
 \nabla_{\partial /\partial t}\, \tilde \xi = \big\{ \tilde \xi, \kk \tilde X - \tilde Y + \sqrt{\sigma}\tilde U\big\},
\end{gather*}
for any $\xi \in {\mathfrak g}$. With the help of the relations~\eqref{relations}, we now obtain
\begin{gather}
 \nabla_{\partial /\partial t}\,\tilde U = \tilde V, \qquad \nabla_{\partial /\partial t}\,\tilde V = -\kk \tilde U,\qquad \nabla_{\partial /\partial t}\,\tilde X = - \tilde H, \nonumber\\
 \nabla_{\partial /\partial t}\,\tilde Y = -\kk \tilde H + \sqrt{\sigma}\tilde V, \qquad \nabla_{\partial /\partial t}\,\tilde H = 2 \kk \tilde X + 2 \tilde Y + \sqrt{\sigma}\tilde U.\label{derivatives}
\end{gather}
With the aid of the Maurer--Cartan equations \eqref{mce}, we may also compute the Lie algebroid bracket of $A(f)$:
\begin{gather}
 [\tilde X,\tilde U] = \tilde H, \qquad [\tilde Y,\tilde U] = \kk \tilde H, \qquad [\tilde H,\tilde U] = -2\kk \tilde X - 2 \tilde Y,\nonumber\\
 [\tilde X,\tilde Y] = \tilde H, \qquad [\tilde H,\tilde X] = 2\tilde X, \qquad [\tilde H,\tilde Y] = -2\tilde Y.\label{brackets}
\end{gather}
One sees, from its definition, that $A(f)^2=\spann\big\{U^2,\tilde X,\tilde H\big\}$, where $U^2=\sqrt{\sigma}\tilde U - \tilde Y$, and we have
\begin{gather*}
 \nabla_{\partial /\partial t}\,U^2= \kk \tilde H,\qquad \nabla_{\partial /\partial t}\,\tilde X = -\tilde H, \qquad \nabla_{\partial /\partial t}\, \tilde H = 2\kk \tilde X + 3\tilde Y + U^2.
\end{gather*}
Evidently, $A(f)^3=\spann\{U^2,\tilde X\}$ and $A(f)^4=\spann{U^4}$, where
\[U^4=U^2+\kk \tilde X = \sqrt{\sigma}\tilde U - \tilde Y +\kk \tilde X.\]
Finally, we compute $\nabla_{\partial /\partial t}\,U^4= \dot \kk \tilde X$, so that $A(f)^5=0$, unless $\kk(t)$ is constant, in which case $A(f)^5=A(f)^4$. Applying Theorem~\ref{sssT}:
\begin{Proposition}\label{gtttP} A path $f \colon I \rightarrow {\mathbb R}^2$ in the equi-affine plane, having constant equi-affine speed $\sigma = s(t)^3 \kappa(t)$, is the orbit of a one-parameter subgroup of area-preserving affine maps precisely when its equi-affine curvature $\kk(t) $ is constant.
\end{Proposition}

In any case, since the transitive Lie algebroid $A^4(f)$ has rank one, we have a canonical identification $A^4(f)\cong T I$, in which the section $U^4$ becomes identified with $\partial /\partial t$, because $\# U^4=\# \sqrt{\sigma}\tilde U =\partial /\partial t$. The canonical representation of $A^4(f)\cong TI$ on $TI$ described in Proposition~\ref{cP} is given by
\begin{gather}
 \bar\nabla_{\partial /\partial t}\partial /\partial t=0.\label{puff}
\end{gather}
This representation equips $I$ with the intrinsic geometry of a~parallelism. Equation \eqref{puff} shows that, with respect to the corresponding parallelism on $\Gamma =f(I)$, the velocity vector~$\dot f$ of the constant equi-affine speed parameterization $f \colon I \rightarrow {\mathbb R}^2 $ is indeed parallel, as claimed earlier.

Now let $Q $ be the symplectic two-form on the radical $\rad {\mathfrak g} \!=\!\spann\{U,V\}$ satisfying \smash{$Q(U,V)\!=\!1$}. Then $Q$ is $\Adjoint$-invariant. In particular, if we define $U^1=\sqrt{\sigma}\tilde U$, with $\tilde U$ as above, so that $\# U^1 = \partial /\partial t$, then
\begin{gather*}
 Q\big(U^1,\nabla_{\partial/\partial t}U^1\big)=\sigma Q\big(\tilde U, \tilde V\big)= \sigma Q(U,V)=\sigma.
\end{gather*}
Next, we let $\mu \colon {\mathfrak g} \rightarrow {\mathbb R}$ denote the quadratic form given by composition of the natural projection ${\mathfrak g} \rightarrow {\mathfrak{sl}}(2,{\mathbb R})$ with the determinant on ${\mathfrak{sl}}(2,{\mathbb R})$. This form is invariant under arbitrary automorphisms of ${\mathfrak g}$. Viewing $\mu $ as a~quadratic form on the trivial bundle ${\mathfrak g} \times I$, we have
\begin{gather*}
 \mu \big( u \tilde U + v \tilde V + a \tilde X + b \tilde Y + c \tilde H \big)=\mu\big( u U + v V + a X + b Y + c H \big) =-c^2-ab.
\end{gather*}
In particular, we have $\mu(U^4)= -\kk$, where $U^4$ is the unique generator of $A(f)^4$ given above such that $\# U^4=\partial /\partial t$. We record for later:
\begin{Lemma}\label{lemmaAA}
 For any $\lambda \ne 0$ there exists an automorphism of ${\mathfrak g} $ preserving $\mu $ and pulling $Q$ back to~$\lambda Q$.
\end{Lemma}
\begin{proof} Put $\epsilon = \sqrt{|\lambda|}$. If $\lambda > 0$, consider the automorphism defined by $U \mapsto \epsilon U$, $V \mapsto \epsilon V$, $X \mapsto X$, $Y \mapsto Y$, $H \mapsto H$. If $\lambda < 0$, take $U \mapsto \epsilon V$, $V \mapsto \epsilon U$, $X \mapsto Y$, $Y \mapsto X$, $H \mapsto -H$.
\end{proof}

Turning our attention to arbitrary infinitesimal immersions:
\begin{Lemma} Let $\omega \colon A \rightarrow {\mathfrak g} $ be any infinitesimal immersion. Then both $A \cap \rad E$ and $A^4$ have rank one and are mapped by the anchor $\# \colon A \rightarrow TI$ onto $TI$. Here $\rad E=\rad {\mathfrak g} \times I$.
\end{Lemma}
\noindent%
\begin{proof} See the proof of Lemma \ref{lemmaB}.
\end{proof}
\begin{Definition} Let $\omega \colon A \rightarrow {\mathfrak g} $ be an infinitesimal immersion of $I$ into ${\mathbb R}^2$ and let $U^1$ and $U^4$ be the unique sections of $A\cap\rad E$ and $A^4$ respectively, such that
 $\# U^1 = \partial / \partial t = \# U^4$. Then the {\df speed} $\sigma \colon I \rightarrow {\mathbb R} $ of $\omega $ is defined by $\sigma =Q\big(U^1, \nabla_{\partial/\partial t}U^1\big)$ and the {\df curvature} $\kk \colon I \rightarrow {\mathbb R} $ of $\omega $ by $\kk=-\mu(U^4)$.
\end{Definition}

Evidently, the speed and curvature are invariants of an infinitesimal immersion. We have just shown that for the logarithmic derivative of a path $f$ with constant equi-affine speed, the speed and curvature coincide with the equi-affine speed and equi-affine curvature of $f$ itself. We are ready to apply Construction Principle~\ref{constructionPrinciple} to prove:
\begin{Theorem} For every smooth function $\kk \colon I \rightarrow {\mathbb R}$ there exists a curve $\Gamma $ in the equi-affine plane, parameterised by some $f \colon I \rightarrow {\mathbb R}^2$ with unit equi-affine speed, whose equi-affine curvature is~$\kk $. The curve is unique up to area-preserving affine transformations.
\end{Theorem}
\begin{proof} Let $\tilde U$, $\tilde V$, $\tilde X$, $\tilde Y$, $\tilde Z$ denote the constant sections of the trivial bundle $E := {\mathbb R}^5 \times I $ and define an algebraic bracket on $E$ by the relations \eqref{relations}. Evidently this makes $E $ into a ${\mathfrak g} $-bundle, with ${\mathfrak g} $ as above. Define a connection $\nabla $ on $E$ by taking $\sigma =1$ in~\eqref{derivatives} and verify that the algebraic bracket is $\nabla $-parallel. Make the subbundle $A \subset E$ spanned by $\tilde U$, $\tilde X$, $\tilde Y$ and $\tilde H$ into a Lie algebroid by defining a bracket by~\eqref{brackets} and an anchor~$\#$ by~\eqref{anchors} (with $\sigma =1$). Then the Maurer--Cartan equations~\eqref{mce} hold. Applying Proposition~\ref{croonP}, we obtain a Lie algebroid morphism $\omega \colon A \rightarrow {\mathfrak g} $ which is readily seen to satisfy the axioms of an infinitesimal immersion.

 We now compute the speed and curvature of $\omega $. Since $Q$, as a~symplectic form on the $\nabla $-invariant bundle $\rad E$, is $\nabla $-parallel, we have
 \begin{gather*}
 \frac{{\rm d}}{{\rm d}t}Q\big(\tilde U,\tilde V\big)=Q\big(\nabla_{\partial /\partial t}\tilde U, \tilde V\big)+Q\big(\tilde U, \nabla_{\partial /\partial t} \tilde V\big)=Q\big(\tilde V, \tilde V\big) - \kk Q\big(\tilde U, \tilde U\big)=0.
 \end{gather*}
Since we have taken $\sigma =1$, the unique section $U^1$ of $A\cap\rad E$ with $\# U^1 = \partial /\partial t$ is $U^1=\tilde U$. The speed of $\omega $ is therefore $Q\big(\tilde U, \nabla_{\partial /\partial t} \tilde U\big)=Q\big(\tilde U, \tilde V\big)$, which is constant, on account of the preceding calculation. By composing $\omega \colon A \rightarrow {\mathfrak g} $ with the appropriate outer automorphism of $\mathfrak g$, we may arrange, by Lemma~\ref{lemmaAA}, for this speed to be one. By a calculation identical to that presented for curves, the unique section $U^4$ of $A^4$ satisfying $\# U^4=\partial /\partial t$ is $U^4 =\tilde U - \tilde Y + \kk \tilde X$. Because the formula for $\mu $ given above must remain the same ($\mu $ is an outer invariant of ${\mathfrak g} $), we compute $\mu\big(U^4\big)=-\kk$. So the curvature of $\omega $ is $\kk$.

By Theorem \ref{infirmT}, $\omega $ has a primitive $f \colon I \rightarrow {\mathbb R}^2$, unique up to area-preserving affine maps. Since this means $\omega $ and $\delta f$ are isomorphic in $\bf{Inf}$, they must have the same invariants. In particular, $f$ must have unit equi-affine speed and the curvature of $\Gamma=f(I)$ must be $\kk$.
\end{proof}

\subsection*{Acknowledgements} We thank Yuri Vyatkin for helpful discussions.

\pdfbookmark[1]{References}{ref}
\LastPageEnding


\begin{thebibliography}{99}
\footnotesize\itemsep=0pt

\bibitem{Armstrong_07}
Armstrong S., Note on pre-{C}ourant algebroid structures for parabolic
 geometries, \href{https://arxiv.org/abs/0709.0919}{arXiv:0709.0919}.

\bibitem{Blaom_06}
Blaom A.D., Geometric structures as deformed infinitesimal symmetries,
 \href{https://doi.org/10.1090/S0002-9947-06-04057-8}{\textit{Trans. Amer. Math. Soc.}} \textbf{358} (2006), 3651--3671,
 \href{https://arxiv.org/abs/math.DG/0404313}{math.DG/0404313}.

\bibitem{Blaom_12}
Blaom A.D., Lie algebroids and {C}artan's method of equivalence, \href{https://doi.org/10.1090/S0002-9947-2012-05441-9}{\textit{Trans.
 Amer. Math. Soc.}} \textbf{364} (2012), 3071--3135, \href{https://arxiv.org/abs/math.DG/0509071}{math.DG/0509071}.

\bibitem{Blaom_13}
Blaom A.D., The infinitesimalization and reconstruction of locally homogeneous
 manifolds, \href{https://doi.org/10.3842/SIGMA.2013.074}{\textit{SIGMA}} \textbf{9} (2013), 074, 19~pages,
 \href{https://arxiv.org/abs/1304.7838}{arXiv:1304.7838}.

\bibitem{Blaom_16b}
Blaom A.D., Cartan connections on {L}ie groupoids and their integrability,
 \href{https://doi.org/10.3842/SIGMA.2016.114}{\textit{SIGMA}} \textbf{12} (2016), 114, 26~pages, \href{https://arxiv.org/abs/1605.04365}{arXiv:1605.04365}.

\bibitem{Blaom_16}
Blaom A.D., Pseudogroups via pseudoactions: unifying local, global, and
 infinitesimal symmetry, \textit{J.~Lie Theory} \textbf{26} (2016), 535--565,
 \href{https://arxiv.org/abs/1410.6981}{arXiv:1410.6981}.

\bibitem{Blaom_F}
Blaom A.D., A characterization of smooth maps into a homogeneous space,
 \href{https://arxiv.org/abs/1702.02717}{arXiv:1702.02717}.

\bibitem{Blaom_H}
Blaom A.D., The Lie algebroid associated with a hypersurface,
 \href{https://arxiv.org/abs/1702.03452}{arXiv:1702.03452}.

\bibitem{Burstall_Calderbank_04}
Burstall F.E., Calderbank D.M.J., Submanifold geometry in generalized flag
 manifolds, \textit{Rend. Circ. Mat. Palermo~(2) Suppl.} \textbf{72} (2004),
 13--41.

\bibitem{Burstall_Calderbank_10}
Burstall F.E., Calderbank D.M.J., Conformal submanifold geometry {I--III},
 \href{https://arxiv.org/abs/1006.5700}{arXiv:1006.5700}.

\bibitem{CannasdaSilva_Weinstein_99}
Cannas~da Silva A., Weinstein A., Geometric models for noncommutative algebras,
 \textit{Berkeley Mathematics Lecture Notes}, Vol.~10, Amer. Math. Soc.,
 Providence, RI, Berkeley Center for Pure and Applied Mathematics, Berkeley,
 CA, 1999.

\bibitem{Cap_Gover_02}
\v{C}ap A., Gover A.R., Tractor calculi for parabolic geometries,
 \href{https://doi.org/10.1090/S0002-9947-01-02909-9}{\textit{Trans. Amer. Math. Soc.}} \textbf{354} (2002), 1511--1548.

\bibitem{Cap_Slovak_09}
\v{C}ap A., Slov\'ak J., Parabolic geometries.~{I}.~Background and general
 theory, \href{https://doi.org/10.1090/surv/154}{\textit{Mathematical Surveys and Monographs}}, Vol.~154, Amer. Math.
 Soc., Providence, RI, 2009.

\bibitem{Crainic_Fernandes_11}
Crainic M., Fernandes R.L., Lectures on integrability of {L}ie brackets, in
 Lectures on {P}oisson Geometry, \textit{Geom. Topol. Monogr.}, Vol.~17, Geom.
 Topol. Publ., Coventry, 2011, 1--107, \href{https://arxiv.org/abs/math.DG/0611259}{math.DG/0611259}.

\bibitem{Dufour_Nguyen_05}
Dufour J.-P., Zung N.T., Poisson structures and their normal forms,
 \href{https://doi.org/10.1007/b137493}{\textit{Progress in Mathematics}}, Vol.~242, Birkh\"auser Verlag, Basel, 2005.

\bibitem{Olver_Fels_98}
Fels M., Olver P.J., Moving coframes. {I}.~{A} practical algorithm,
 \href{https://doi.org/10.1023/A:1005878210297}{\textit{Acta Appl. Math.}} \textbf{51} (1998), 161--213.

\bibitem{Olver_Fels_99}
Fels M., Olver P.J., Moving coframes. {II}.~{R}egularization and theoretical
 foundations, \href{https://doi.org/10.1023/A:1006195823000}{\textit{Acta Appl. Math.}} \textbf{55} (1999), 127--208.

\bibitem{Mackenzie_05}
Mackenzie K.C.H., General theory of {L}ie groupoids and {L}ie algebroids,
 \href{https://doi.org/10.1017/CBO9781107325883}{\textit{London Mathematical Society Lecture Note Series}}, Vol.~213, Cambridge
 University Press, Cambridge, 2005.

\bibitem{Olver_05}
Olver P.J., A survey of moving frames, in Computer Algebra and Geometric
 Algebra with Applications (6th International Workshop, IWMM 2004, Shanghai,
 China, May 19--21, 2004 and International Workshop, GIAE 2004, Xian, China,
 May 24--28, 2004), \href{https://doi.org/10.1007/11499251_11}{\textit{Geom. Topol. Monogr.}}, Vol.~17, Editors H.~Li,
 P.J.~Olver, G.~Sommer, Springer, Berlin~-- Heidelberg, 2005, 105--138.

\bibitem{Sharpe_97}
Sharpe R.W., Differential geometry. {C}artan's generalization of {K}lein's
 {E}rlangen program, \textit{Graduate Texts in Mathematics}, Vol.~166,
 Springer-Verlag, New York, 1997.

\bibitem{XuXiaomeng_14}
Xu X., Twisted {C}ourant algebroids and coisotropic {C}artan geometries,
 \href{https://doi.org/10.1016/j.geomphys.2014.03.002}{\textit{J.~Geom. Phys.}} \textbf{82} (2014), 124--131, \href{https://arxiv.org/abs/1206.2282}{arXiv:1206.2282}.

\end{thebibliography}
\end{document}